\documentclass[11pt]{amsart}
\usepackage{amsmath, amsthm, amsfonts, graphicx, ifpdf}
\usepackage{mathrsfs,epsfig}
\usepackage[dvipsnames,usenames]{color}
\usepackage[colorlinks=true,
            linkcolor=blue,
            citecolor=blue,
            urlcolor=blue,
            bookmarksnumbered=true]{hyperref}
\hypersetup{
   unicode=false,          
   pdftoolbar=true,        
bookmarksnumbered=true,
   pdfmenubar=true,        
   pdffitwindow=false,     
   pdfstartview={},    
   pdftitle={},    
   pdfauthor={},     
   pdfsubject={},   
   pdfcreator={},   
   pdfproducer={}, 
   pdfkeywords={}, 
   pdfnewwindow=true,      
   colorlinks=true,       
   linkcolor=blue,          
   citecolor=blue,        
   filecolor=blue,      
   urlcolor=blue          
}

\usepackage[T1]{fontenc}
\usepackage{lmodern,mathtools,amssymb,microtype,needspace}
\usepackage[a4paper,textwidth=158mm,textheight=238mm,centering]{geometry}
\numberwithin{equation}{section}
\newtheorem{mainthm}{Theorem}
\newtheorem{theorem}{Theorem}[section]
\newtheorem{lemma}[theorem]{Lemma}
\newtheorem{proposition}[theorem]{Proposition}
\newtheorem{corollary}[theorem]{Corollary}
\theoremstyle{remark}
\newtheorem*{remark}{Remark}
\newtheorem{ques}[theorem]{Question}

\newtheorem*{rem*}{Remark}

\newcommand{\ke}{ds^2_{\mathrm{KE}}}

\newcommand{\be}{\begin{equation}}
\newcommand{\ene}{\end{equation}}
\newcommand{\br}{\begin{remark}}
\newcommand{\er}{\end{remark}}
\newcommand{\bl}{\begin{lem}}
\newcommand{\el}{\end{lem}}
\newcommand{\bcor}{\begin{cor}}
\newcommand{\ecor}{\end{cor}}
\newcommand{\bpro}{\begin{pro}}
\newcommand{\epro}{\end{pro}}
\newcommand{\ben}{\begin{enumerate}}
\newcommand{\een}{\end{enumerate}}
\newcommand{\bp}{\begin{proof}}
\newcommand{\ep}{\end{proof}}
\newcommand{\bpo}{\begin{pro}}
\newcommand{\epo}{\end{pro}}
\newcommand{\beq}{\begin{equation*}}
\newcommand{\eeq}{\end{equation*}}
\newcommand{\bear}{\begin{eqnarray}}
\newcommand{\eear}{\end{eqnarray}}
\newcommand{\beqar}{\begin{eqnarray*}}
\newcommand{\eeqar}{\end{eqnarray*}}
\newcommand{\brem}{\begin{remark*}}
\newcommand{\erem}{\end{remark*}}
\newcommand{\bt}{\begin{theorem}}
\newcommand{\et}{\end{theorem}}

\newcommand{\C}{\mathbb C}
\newcommand{\R}{\mathbb R}
\newcommand{\Q}{\mathbb Q}
\newcommand{\D}{\mathbb D}
\newcommand{\B}{\mathbb B}
\newcommand{\bnum}[1]{b^{(2)}_{#1}}
\newcommand{\Harm}{\mathcal H}
\newcommand{\Nalg}{\mathcal N}
\newcommand{\quot}{\backslash}
\newcommand{\frg}{\mathfrak g}
\newcommand{\frk}{\mathfrak k}
\newcommand{\frp}{\mathfrak p}
\newcommand{\frt}{\mathfrak t}
\newcommand{\ip}[2]{\langle #1,#2\rangle}
\DeclareMathOperator{\Aut}{Aut}
\DeclareMathOperator{\Vol}{Vol}
\DeclareMathOperator{\Ric}{Ric}

\DeclareMathOperator{\rank}{rank}

\DeclareMathOperator{\covol}{covol}
\DeclareMathOperator{\Int}{Int}
\DeclareMathOperator{\Span}{span}
\DeclareMathOperator{\dist}{dist}
\DeclareMathOperator{\CE}{CE}
\DeclareMathOperator{\id}{id}
\newcommand{\sq}{\mathsf{s}}

\numberwithin{equation}{section}

\allowdisplaybreaks

\def\XXint#1#2#3{{\setbox0=\hbox{$#1{#2#3}{\int}$}
    \vcenter{\hbox{$#2#3$}}\kern-.5\wd0}}

\makeatletter
\def\@citestyle{\m@th\upshape\mdseries}
\def\citeform#1{{\bfseries#1}}
\def\@cite#1#2{{%
  \@citestyle[\citeform{#1}\if@tempswa, #2\fi]}}
\@ifundefined{cite }{%
  \expandafter\let\csname cite \endcsname\cite
  \edef\cite{\@nx\protect\@xp\@nx\csname cite \endcsname}%
}{}
\makeatother

\title[Rigidity of bounded domains]{Rigidity of complex bounded domains with finite-volume quotients}

\author{Le Kuai}
\address{Qiuzhen College, Tsinghua University, Beijing, 100084, P. R. China}
\email{kuail23@mails.tsinghua.edu.cn}

\author{Kefeng Liu}
\address{Mathematical Sciences Research Center, Chongqing University of Technology,
 Chongqing 400054, P. R. China;\\
 Department of Mathematics, University of California at Los Angeles, Los Angeles,
 CA 90095-1555, USA}
\email{liu@math.ucla.edu}

\author{Yunhui Wu}
\address{Department of Mathematical Sciences \& The Yau Mathematical Sciences Center, 
Tsinghua University, Beijing 100084, P. R. China}
\email{yunhui\_wu@mail.tsinghua.edu.cn}

\date{}
\begin{document}
\begin{abstract} We prove that any bounded convex domain in $\C^n$ admitting a finite-volume quotient is biholomorphic to a bounded symmetric domain. We first show that such a convex domain always has a nontrivial one-parameter family of automorphisms; the proof then relies on our solution to Kazhdan's conjecture for contractible HHR/USq domains with finite-volume quotients. More precisely, we prove that any contractible bounded domain that has the uniform squeezing property and admits a finite-volume quotient is biholomorphic to a product of a bounded symmetric domain and a contractible complex manifold with discrete automorphism group.  The quotient groups are discrete and may have finite stabilizers; volume is measured with respect to the complete K\"ahler--Einstein metric. As a further application, we also show that the automorphism groups of Bers--Griffiths domains in $\C^2$, which are not biholomorphic to the bidisk, are discrete.   
\end{abstract}
\maketitle

\section{Introduction}

In 1989, Frankel~\cite{Fra89} proved that a bounded convex domain in $\C^n$ with a compact quotient is biholomorphic to a \emph{bounded symmetric domain}, addressing a conjecture of Yau~\cite{Yau87}. The finite-volume analogue has remained a folklore conjecture for over three decades; see \emph{e.g.} Siu's survey~\cite[p.~124]{Siu91} and \cite[p.~163]{Zim21}. Our first theorem settles this conjecture.
\begin{mainthm}\label{thm:convex}
A bounded convex domain in $\C^n$ with a finite-volume quotient
is biholomorphic to a bounded symmetric domain.
\end{mainthm}

Throughout the paper, volume is measured with respect to the complete
K\"ahler--Einstein metric, normalized by $\Ric(\ke)=-\ke$. By a
\emph{finite-volume quotient} of a domain $U \subset \C^n$, we mean a quotient
$\Gamma\quot U$ by a discrete subgroup $\Gamma<\Aut(U)$ with finite
orbifold volume.  Our main interest is in noncompact quotients; the arguments also apply
to compact quotients. The action is properly discontinuous, and its point
stabilizers are finite.

\begin{rem*}
For $3g+n\ge5$, the Teichm\"uller space $\mathcal T_{g,n}$ of
Riemann surfaces of genus $g$ with $n$ punctures is not biholomorphic
to a bounded convex domain. Indeed, its moduli-space quotient has
finite K\"ahler--Einstein volume~\cite{McM00, LSY04}, whereas the
automorphism theorems of Royden and Earle--Kra~\cite{Roy71,EK74}
imply that $\Aut(\mathcal T_{g,n})$ is discrete in this range.
Theorem~\ref{thm:convex} would turn any bounded convex realization
into a bounded symmetric domain, contradicting this discreteness.
This recovers Siu's conjecture~\cite[Conjecture~3.7]{Siu91}.
Using deep Teichm\"uller theory, Markovic~\cite{Mar18} first proved the result for closed surfaces;
the extension to all $\mathcal T_{g,n}$ of complex dimension at least
two is given in \cite[Theorems~1.6--1.7]{GM20}. Theorem~\ref{thm:convex}
gives a proof from finite-volume rigidity and the classical
automorphism theorem.
\end{rem*}

The proof of Theorem~\ref{thm:convex} combines the construction of a
nontrivial one-parameter subgroup of holomorphic automorphisms with
a product decomposition for a larger class of domains. The second
part is related to Kazhdan's conjecture: an irreducible bounded domain
with a compact quotient and a nontrivial one-parameter group of
holomorphic automorphisms should be biholomorphic to a bounded
symmetric domain. Here \emph{irreducible} means that the domain is not
biholomorphic to a product of two positive-dimensional complex
manifolds. Frankel's convex-domain theorem was an early case of the
conjecture~\cite{Fra89}. Nadel~\cite{Nad90} established the required semisimplicity of the identity component of the automorphism group, and
Frankel~\cite{Fra95} completed the compact case by harmonic-map
methods. The extension of Kazhdan's conjecture to finite-volume quotients is a well-known \emph{open} problem; see \emph{e.g.}~\cite{Siu91, Fra95}. One may also see the work of Yau--Zheng~\cite{YZ93} for related questions. We establish the finite-volume extension for contractible
domains with the following uniformity property.

Following Liu--Sun--Yau~\cite{LSY04,LSY05}, a complex manifold $U$
of dimension $n$ is \emph{holomorphic homogeneous regular} (HHR) if
there is a constant $a\in(0,1]$ such that, for every $p\in U$, a
biholomorphic embedding $f_p:U\to\C^n$ satisfies
\[
              f_p(p)=0,\qquad B(0;a)\subset f_p(U)\subset B(0;1),
\]
where $B(0;r)$ is the Euclidean ball of radius $r$ in $\C^n$.
Yeung~\cite{Yeu09} calls this the \emph{uniform squeezing property}
(USq); we use the abbreviation HHR/USq. Its motivation goes back to
the uniform regularity conditions in Morrey's work on higher-dimensional
Plateau problems~\cite[Chapter~10]{Mor08}. Examples include bounded homogeneous domains and bounded domains
with compact manifold quotients~\cite{Yeu09}; Teichm\"uller spaces, through the uniform estimates for the Bers embeddings~\cite{Ber60}; and bounded convex domains~\cite{Fra91,KZ16}. The class is closed
under finite products. Every HHR/USq manifold has a bounded pseudoconvex
realization~\cite{Yeu09}, so the existence theorems of
Cheng--Yau~\cite{CY80} and Mok--Yau~\cite{MY83} give a complete
K\"ahler--Einstein metric. Our second main result is as follows.
\begin{mainthm}\label{thm:hhr}
Let $U$ be a contractible bounded $\mathrm{HHR/USq}$ domain with a finite-volume
quotient. Then
\[
                              U\cong S\times Y
\]
biholomorphically, where $S$ is a bounded symmetric domain, possibly
a point, and $Y$ is a contractible complex manifold with discrete
automorphism group.
\end{mainthm}

If $U$ is irreducible and admits a nontrivial one-parameter group of
automorphisms, then $S$ has positive dimension and $Y$ must be a point.
Thus Theorem~\ref{thm:hhr} gives the finite-volume form of Kazhdan's
conjecture in the contractible HHR/USq class.

For a bounded convex domain, the factor $Y$ admits a bounded convex
realization~\cite[Corollary~5.5]{Fra89}. Our splitting argument also
gives a finite-volume quotient of this factor. We prove that every
positive-dimensional bounded convex domain with such a quotient
admits a nontrivial one-parameter group of holomorphic automorphisms.
Applied to $Y$, this contradicts the discreteness of $\Aut(Y)$ unless
$Y$ is a point, and proves Theorem~\ref{thm:convex}.

We also apply Theorem~\ref{thm:hhr} to Bers--Griffiths domains.
Griffiths~\cite{Gri71} constructed contractible bounded domains in
arbitrary dimensions as universal covers of suitable Zariski open
subsets of algebraic varieties. In complex dimension two, his
construction gives a smooth algebraic family $\pi:M\to C$ of
finite-type hyperbolic curves over a finite-type hyperbolic curve.
The precise family hypothesis is recalled in
Section~\ref{subsec:griffiths-families}. Bers' relative uniformization
realizes $D=\widetilde M$ as a bounded pseudoconvex domain with a
disk fibration over $\D$. We call these two-dimensional universal
covers \emph{Bers--Griffiths domains}; related uniformization and
automorphism questions are studied in~\cite{GD08,GR18,Ima83,Sha77}. Fornaess and Kim~\cite[Problem~18]{FK15} recorded the corresponding
automorphism question. With the bidisk excluded, it takes the form:

\begin{ques}\label{op-2}
Let $D$ be a two-dimensional Bers--Griffiths domain that is not biholomorphic to the
bidisk. Is $\Aut(D)$ discrete?
\end{ques}

Shabat~\cite[Theorem~3]{Sha77} proved discreteness, apart from the
bidisk, when the base or the fibers are compact. The difficult part of Question~\ref{op-2} is the case that both the base and the fibers of $M$ are open surfaces. We give a positive answer to Question~\ref{op-2} through
Theorem~\ref{thm:hhr}, after establishing the HHR/USq property for
these two-dimensional Bers--Griffiths domains.

\begin{mainthm}\label{mt-3}
Let $D$ be a two-dimensional Bers--Griffiths domain. If $D$ is not
biholomorphic to the bidisk $\D\times\D$, then $\Aut(D)$ is discrete.
\end{mainthm}

The additional input is quantitative: for every $p\in D$ there is a
holomorphic embedding $\Psi_p:D\hookrightarrow\B^2$ with
\[
                    \Psi_p(p)=0,\qquad
                    \B_{1/5}^2\subset\Psi_p(D).
\]
We prove this for all smooth families of the above type in
Theorem~\ref{thm:griffiths-squeezing}; it answers the surface case of
Zimmer's question~\cite[Question~1.18]{Zim26}. For a Bers--Griffiths domain, the algebraic
quotient $M$ has finite K\"ahler--Einstein volume by
\cite[Proposition~7.3]{Gri71}. Theorem~\ref{thm:hhr} then reduces
Theorem~\ref{mt-3} to the ball and the bidisk. The first is excluded
by Imayoshi--Nishimura~\cite[Theorem~1]{IN05}, and the second by
hypothesis.

\subsection*{Methods}
We now describe the two parts of the proof of Theorem~\ref{thm:hhr}. Write
$G=\Aut_0(U)$ and $\Gamma_0=\Gamma\cap G$. The passage from finite
volume to the structure of $G$ is as follows:
the frame-bundle argument makes $\Gamma_0$ a lattice in $G$, and
Prasad's rank formula, together with nonvanishing of group $L^2$
cohomology, controls the solvable radical and the center. Using a
conjugation argument, we obtain
\[
                         \Gamma_b=\Gamma_0\times\Lambda
\]
after passage to a finite-index subgroup. We give these arguments with
finite stabilizers included. In particular, compact normal subgroups
are excluded by a bounded plurisubharmonic function on the quotient.

The next step is to relate this group decomposition to the geometry of
$U$. Let $K<G$ be a maximal compact subgroup. We replace $\Gamma_0$
by a Cartan-invariant uniform lattice $L$, leaving $\Lambda$ unchanged.
An estimate for the first jet of a Killing field turns the adjoint
spectral gap of $L$ into a uniform displacement bound on $U$. It follows
that the translates of the Cartan fixed manifold are uniformly
separated. A tube of fixed radius therefore descends to the quotient,
and integration along its normal fibers proves that the fixed manifold
has finite induced quotient volume. Torus localization and the
Cheeger--Gromov comparison theorem now give
\[
                         \dim_\R U=\dim_\R(G/K)+\dim_\R U^K.
\]
The orbit map $(G/K)\times U^K\to U$ is consequently a global
biholomorphism. 

For Theorem~\ref{thm:convex}, by Theorem~\ref{thm:hhr} it remains to show that a
positive-dimensional bounded convex domain with a finite-volume
quotient has a nontrivial one-parameter group of automorphisms. We find
a complex extreme boundary point at which the span of the normal cone
is totally real. The existence of this point follows from the positivity
of an integral over the unit normal bundle. A global peak function and
finite volume put this point in an orbit limit set. Supporting
hyperplanes then control rescaled translations, both on compact sets
and along their trajectories, and give a nonzero complete holomorphic
vector field. Applying this result to the convex realization of $Y$
shows that $Y$ is a point.

For Theorem~\ref{mt-3} about Bers--Griffiths domains,  by Theorem~\ref{thm:hhr} it remains to show that such domains are HHR/USq. The reference fiber in Bers' construction is
chosen anew at each point. Principal normalization gives a uniform
outer bound for every fiber. Schwarz's lemma in the base variable
then supplies a uniform neighborhood of the origin inside the total
domain. A fixed linear change of coordinates yields the squeezing
constant $1/5$. This argument uses that both the base and the fibers
are one-dimensional.

\subsection*{Plan of the paper}
Section~\ref{sec:prelim} contains the analytic preliminaries.
Section~\ref{sec:structure} establishes the group structure, and
Section~\ref{sec:fixed-volume} proves finite volume for the fixed
manifold. Theorem~\ref{thm:hhr} is proved in
Section~\ref{sec:splitting}. Sections~\ref{sec:convex-geometry}
and~\ref{sec:flows}  develop the convex
boundary argument and prove Theorem~\ref{thm:convex}.  Section~\ref{sec:griffiths} establishes
the uniform squeezing property for Bers--Griffiths domains and deduces
Theorem~\ref{mt-3} from Theorem~\ref{thm:hhr}.

\subsection*{AI disclosure} All mathematical ideas were developed by the authors. AI tools were also used to assist in the preparation and editing of this revised manuscript, including simplifying the presentation of proofs, completing references, and polishing the writing. The authors take full responsibility for the final manuscript.

\subsection*{Acknowledgments}
The second and third named authors thank Professors W. Ballmann, S. Krantz, B. Wong, and S. T. Yau for their interest, invaluable comments, and suggestions on the first version of this paper, which was written eight years ago. The second named author is partially supported by the National Key Research and Development Program of China (NKPs) [Grant Number 2024YFA1013201], the Natural Science Foundation of Chongqing (NSFCQ) [Grant Number CSTB2024NSCQ-LZX0040], and the Special Project of Chongqing Municipal Science and Technology Bureau [Grant Number 2025CCZ015].
The third named author is partially supported by the National Key R \& D Program of China (2025YFA1017500) and NSFC grants No. 12361141813 and 12425107. 

\section{Notation and preliminaries}\label{sec:prelim}

Until the end of Section~\ref{sec:splitting}, $U\subset\C^n$ is a
contractible bounded HHR/USq domain, $n\ge2$, and
$\Gamma<\Aut(U)$ is a discrete subgroup such that
\[
              M=\Gamma\quot U,
              \qquad \Vol_{\mathrm{KE}}(M)<\infty.
\]
Thus $M$ is understood as a Riemannian orbifold. We write
\[
                      G=\Aut_0(U),\qquad \Gamma_0=\Gamma\cap G.
\]
The automorphism group is given the compact-open topology. We use
$d_{\mathrm{KE}}$ for the distance induced by $\ke$ and
$d_U^{\mathrm K}$ for the Kobayashi distance on $U$; its infinitesimal
metric is denoted by $k_U$. Euclidean balls are written
$\B_r^j\subset\C^j$, with $\B^j=\B_1^j$ and
$\D_r=\B_r^1$, $\D=\D_1$. The factors in Theorem~\ref{thm:hhr} are denoted by $S$ and $Y$,
with points $(s,y)$; the letter $f$ is used for holomorphic maps.
We write $\mathcal X$ for the auxiliary manifold in the dimension
comparison and $C$ for the base curve in Section~\ref{sec:griffiths}.
In the convex boundary argument, $\mathcal K$ denotes a convex body
and $V_p$ the real span of its normal cone at $p$; $K$ is reserved
for compact subgroups of the automorphism group.

\subsection{The K\"ahler--Einstein metric and quotient volume}

A bounded HHR/USq domain is pseudoconvex~\cite[Theorem~1]{Yeu09} and admits a complete K\"ahler--Einstein metric~\cite{CY80, MY83}. We use the
normalization
\[
                         \Ric(\ke)=-\ke,
\]
for which the metric is unique. After rescaling Yeung's
normalization, the curvature tensor of $\ke$ and all its covariant
derivatives are uniformly bounded, and the injectivity radius of
$(U,\ke)$ is bounded below by a positive constant
\cite[Theorem~2(b),(d)]{Yeu09}. One may also see \cite{LSY04, LSY05}. We refer to these properties as
bounded geometry.

The group $\Aut(U)$ is a second-countable Lie group acting properly
on $U$~\cite{Kob98}. By uniqueness of $\ke$, this action is isometric.
Thus $\Gamma$ is countable and has finite point stabilizers. Let
$q:U\to\Gamma\quot U$ be the quotient map. Since $(U,d_{\mathrm{KE}})$
is proper and every $\Gamma$-orbit is closed, for all $x,y\in U$ we have
\[
d_{\Gamma\quot U}\bigl(q(x),q(y)\bigr)
   =\min_{\gamma\in\Gamma}d_{\mathrm{KE}}(x,\gamma y).
\]
It follows that
\[
\overline B_{\Gamma\quot U}\bigl(q(x),R\bigr)
   =q\bigl(\overline B_U(x,R)\bigr).
\]
The ball on the right is compact by Hopf--Rinow. Hence the quotient
metric is proper, and therefore complete.

For every finite-index subgroup $\Gamma'<\Gamma$, integration in
finite quotient charts gives
\begin{equation}\label{eq:index}
       \Vol(\Gamma'\quot U)=[\Gamma:\Gamma']\Vol(\Gamma\quot U).
\end{equation}
One may also compute these integrals on the free locus. The fixed set of
a nonidentity holomorphic automorphism is a proper analytic subset.
The union of these sets over the countable group $\Gamma$ has measure
zero, so the action is free almost everywhere.

A torsion-free discrete subgroup acts freely, since its stabilizers are
finite. Conversely, suppose that a discrete subgroup contains a
nonidentity element of finite order. A suitable nontrivial power has
prime order $p$. Since the contractible manifold $U$ is mod-$p$
acyclic, Smith theory gives a fixed point for this power~\cite{Bre72}.
Thus a discrete action on $U$ is free if and only if the acting group
is torsion-free. We will use this observation only for an auxiliary
lattice; $\Gamma$ itself may have torsion.

\subsection{A finite-volume Liouville lemma}

The following form of the Liouville argument will be used for both the
automorphism group and boundary orbit limits.

\begin{lemma}\label{lem:liouville}
A bounded nonnegative locally Lipschitz weakly subharmonic function
on a connected complete finite-volume Riemannian orbifold is constant.
\end{lemma}

\begin{proof}
Let $M$ be the orbifold and $u$ the given function.
By completeness, we may choose compactly supported Lipschitz
cutoffs $\chi_R$ such that $\chi_R=1$ on $B_R$, $\chi_R=0$
outside $B_{2R}$, and
\[
0\le\chi_R\le1,\qquad |\nabla\chi_R|\le C/R
\quad\text{almost everywhere}.
\]
Testing $\Delta u\ge0$ against the nonnegative function
$\chi_R^2u$ gives
\[
\begin{aligned}
0
&\le -\int_M
   \bigl\langle\nabla u,\nabla(\chi_R^2u)\bigr\rangle\\
&=-\int_M\chi_R^2|\nabla u|^2
  -2\int_M\chi_Ru
       \langle\nabla\chi_R,\nabla u\rangle .
\end{aligned}
\]
Consequently, Young's inequality yields
\begin{equation}\label{eq:energy}
\int_M\chi_R^2|\nabla u|^2
\le 4\int_Mu^2|\nabla\chi_R|^2
\le
\frac{4C^2\|u\|_\infty^2}{R^2}\Vol(M)
\longrightarrow0.
\end{equation}

The test is justified by smooth approximations in
quotient charts, averaging over each finite chart group, and
an orbifold partition of unity. In a chart $\widetilde V/H$,
the notation $|H|$ denotes the order of the finite group $H$:
orbifold integration is defined as $|H|^{-1}$ times integration
on $\widetilde V$. Thus, integration by parts carries the same
factor on both sides. For every compact set $Q\subset M$, we have
$\chi_R=1$ on $Q$ for all sufficiently large $R$.
Equation~\eqref{eq:energy}
therefore implies $\nabla u=0$ almost everywhere on $M$.
Since $u$ is locally Lipschitz and $M$ is connected, $u$ is
constant.
\end{proof}

\begin{lemma}\label{lem:sup}
Let $\{f_j:V\to\C^m\}$ be a uniformly bounded family of holomorphic
maps on a complex manifold $V$. Then $u=\sup\limits_j\|f_j\|^2$ is bounded,
locally Lipschitz, and plurisubharmonic.
\end{lemma}

\begin{proof}
This should be well-known to experts. We outline a proof here. Choose $A<\infty$ such that
$\|f_j(z)\|\le A$ for all $j$ and $z\in V$. Fix coordinate balls
$B'\Subset B\Subset V$. Cauchy estimates give a constant $C$,
independent of $j$, such that $\|Df_j\|\le C$ on $B'$.
Since $B'$ is convex in these coordinates,
\[
\begin{aligned}
\bigl|\|f_j(z)\|^2-\|f_j(w)\|^2\bigr|
&\le
\bigl(\|f_j(z)\|+\|f_j(w)\|\bigr)
\|f_j(z)-f_j(w)\|\\
&\le 2AC|z-w|,
\qquad z,w\in B'.
\end{aligned}
\]
Taking suprema gives $|u(z)-u(w)|\le2AC|z-w|$.
Thus $u$ is locally Lipschitz and $0\le u\le A^2$.
Each $\|f_j\|^2$ is plurisubharmonic. The upper-envelope
theorem~\cite[Chapter~I, Theorem~5.7]{Dem12} makes the upper
semicontinuous regularization of $u$ plurisubharmonic; continuity
identifies this regularization with $u$.
\end{proof}

On a K\"ahler manifold, a plurisubharmonic function is weakly
subharmonic. Hence a function as in Lemma~\ref{lem:sup} is constant
if it descends to a complete finite-volume quotient.


\subsection{\texorpdfstring{$L^2$}{L2} cohomology and dimension comparison}

We write $\bnum{k}(\Pi)$ for the group $L^2$-Betti numbers,
$\Nalg(\Pi)$ for the group von Neumann algebra with its
canonical trace, and $\Harm^k_{(2)}(\mathcal X)$ for the space of
square-integrable harmonic $k$-forms~\cite{Luc02}.

\begin{proposition}\label{prop:l2-dimension}
Let $U\subset\C^n$ be a contractible bounded $\mathrm{HHR/USq}$ domain,
and let $\Pi<\Aut(U)$ be discrete with
$\Vol_{\mathrm{KE}}(\Pi\quot U)<\infty$.
Then the following statements hold.
\begin{enumerate}
\item
The middle $L^2$-Betti number satisfies
\begin{equation}\label{eq:l2-nonvanishing}
                     0<\bnum{n}(\Pi)<\infty.
\end{equation}
Every amenable normal subgroup of $\Pi$ is finite.

\item
Suppose that $\Pi$ also acts properly by orientation-preserving
isometries on a complete oriented Riemannian manifold $\mathcal X$
without boundary, of bounded geometry, with
\[
\widetilde H_*(\mathcal X;\Q)=0,
\qquad
\Vol(\Pi\quot \mathcal X)<\infty.
\]
Then $\dim_\R \mathcal X=2n$. No freeness assumption is imposed.
\end{enumerate}
\end{proposition}

\begin{proof}
By~\cite[Theorem~2(b),(d)]{Yeu09}, the complete metric $\ke$ has
bounded geometry. The holomorphic volume form
\[
                    \Omega_U=dz_1\wedge\cdots\wedge dz_n
\]
has finite $L^2$-norm, because $|\Omega_U|_{\mathrm{KE}}^2
d\operatorname{vol}_{\mathrm{KE}}$ is a fixed positive multiple of
Euclidean volume and $U$ is bounded. A holomorphic top form is
closed, and the complex-linear Hodge star of a K\"ahler metric
satisfies $*\Omega_U=c_n\Omega_U$, where $|c_n|=1$. Thus
$d\Omega_U=d^*\Omega_U=0$. Completeness identifies the minimal and
maximal de Rham extensions, so
$0\ne\Omega_U\in\Harm^n_{(2)}(U)$.

We now pass from equivariant analytic cohomology to the group
invariant. The action of $\Pi$ on $U$ is proper and has finite
stabilizers. Choose a locally finite $\Pi$-equivariant triangulation
and barycentrically subdivide it so that every cell stabilizer fixes
the cell pointwise. If a cell has stabilizer $H$, its orbit contributes
the cellular module
\[
             \C[\Pi/H]\cong \C\Pi\otimes_{\C H}\C.
\]
This module is projective over $\C\Pi$, because $H$ is finite and the
trivial $\C H$-module is projective by averaging over $H$. Since $U$
is contractible, the augmented cellular complex
\[
                    C_*(U;\C)\longrightarrow\C
\]
is therefore a projective resolution of the trivial $\C\Pi$-module.
Consequently, in the algebraic definition of equivariant
$L^2$-Betti numbers~\cite[Chapter~6]{Luc02},
\[
 \dim_{\Nalg(\Pi)}
 H_k\!\left(\Nalg(\Pi)\otimes_{\C\Pi}C_*(U;\C)\right)
 =\bnum{k}(\Pi).
\]
Here the tensor product is algebraic; no finite-type assumption
or $\ell^2$ completion of the cellular resolution is used.
The Cheeger--Gromov analytic comparison theorem, which includes
proper nonfree discrete actions~\cite[Theorem~5.1]{CG86}, now gives
\begin{equation}\label{eq:comparison}
\bnum{k}(\Pi)
=
\dim_{\Nalg(\Pi)}\Harm^k_{(2)}(U)
<\infty.
\end{equation}
To justify positivity, the $\Pi$-action is free off a null set,
as observed above. A measurable fundamental domain therefore
identifies the space of $L^2$ forms with a Hilbert space of the
form $\ell^2(\Pi)\mathbin{\widehat\otimes}\mathscr V$ for a Hilbert space $\mathscr V$.
The harmonic forms constitute a closed $\Pi$-invariant subspace
containing $\Omega_U\ne0$. Its orthogonal projection is therefore
nonzero, and the induced faithful semifinite trace gives this
subspace positive dimension, proving~\eqref{eq:l2-nonvanishing}.
An infinite amenable normal subgroup would force all
$\bnum{k}(\Pi)$ to vanish~\cite{CG86,Luc02}.
This proves the first assertion.

For the second assertion,
$\widetilde H_*(\mathcal X;\Q)=0$ implies $\widetilde H_*(\mathcal X;\C)=0$.
Properness again makes the cell stabilizers finite. Thus the same
cellular argument gives a projective resolution over $\C\Pi$,
and the equivariant topological $L^2$-Betti numbers of $\mathcal X$ are
those of $\Pi$. Contractibility of $\mathcal X$ is not needed.

To make the analytic comparison explicit when the action has a
kernel, let
\[
                    A=\ker(\Pi\curvearrowright \mathcal X).
\]
Properness makes $A$ finite, and $\overline\Pi=\Pi/A$ is a discrete
subgroup of $\operatorname{Isom}(\mathcal X)$. Counting the full stabilizers
in orbifold volume gives
$\Vol(\overline\Pi\quot \mathcal X)=|A|\Vol(\Pi\quot \mathcal X)<\infty$.
Apply the cellular argument and Cheeger--Gromov comparison to
this effective action. Since $A$ acts trivially on harmonic
forms, the finite-kernel formulas~\cite[Chapter~6]{Luc02} give
\[
\begin{aligned}
 \bnum{k}(\Pi)
 &=\frac1{|A|}\bnum{k}(\overline\Pi)
  =\frac1{|A|}\dim_{\Nalg(\overline\Pi)}\Harm^k_{(2)}(\mathcal X)\\
 &=\dim_{\Nalg(\Pi)}\Harm^k_{(2)}(\mathcal X).
\end{aligned}
\]

Set $d=\dim_\R \mathcal X$ and
\[
                 \mathcal I=\{k\in\mathbf Z_{\ge0}:\bnum{k}(\Pi)>0\}.
\]
The set $\mathcal I$ is nonempty by~\eqref{eq:l2-nonvanishing}. It is finite:
the comparison with $U$ gives $\bnum{k}(\Pi)=0$ for $k>2n$.

The action of $\Pi$ on $U$ is holomorphic and hence
orientation-preserving. Hodge star for $\ke$ is therefore a unitary
$\Pi$-equivariant isomorphism, and~\eqref{eq:comparison} gives
\begin{equation}\label{eq:duality-u}
                 \bnum{k}(\Pi)=\bnum{2n-k}(\Pi)
                 \qquad(0\le k\le2n).
\end{equation}
The same argument on the oriented $d$-manifold $\mathcal X$ gives
\begin{equation}\label{eq:duality-aux}
                 \bnum{k}(\Pi)=\bnum{d-k}(\Pi)
                 \qquad(0\le k\le d).
\end{equation}
In particular, comparison with $\mathcal X$ shows that $\mathcal I\subset[0,d]$.
Thus the two reflections
\[
                  r_U(k)=2n-k,
                  \qquad
                  r_{\mathcal X}(k)=d-k
\]
both permute the nonempty finite set $\mathcal I$. Their composition acts on
$\mathcal I$ as the translation
\[
                  r_U\circ r_{\mathcal X}(k)=k+(2n-d).
\]
A nonempty finite subset of $\mathbf Z$ cannot be invariant under a
nonzero translation. Hence $2n-d=0$, as claimed.
\end{proof}

\begin{rem*}
Only the nonvanishing in~\eqref{eq:l2-nonvanishing} is used above.
In particular, the argument neither asserts that the
K\"ahler--Einstein form has a bounded primitive nor requires
$L^2$-cohomology to vanish away from the middle degree. Uniform
equivalence of canonical metrics alone would not provide such a
primitive; the dimension conclusion instead follows from the two
Hodge-star symmetries~\eqref{eq:duality-u} and~\eqref{eq:duality-aux}.
\end{rem*}


\section{Semisimple, no compact factors, and centerless}\label{sec:structure}

In this section, we prove that $G=\Aut_0(U)$ is centerless and semisimple
without compact factors, and that $\Gamma$ splits after passage to
a subgroup of finite index. 

\subsection{The lattice \texorpdfstring{$\Gamma_0$}{Gamma0}}

We begin with a lattice argument. It uses only the invariant metric and
finite quotient volume.
\begin{lemma}\label{lem:lattice} 
Let $U$ be a bounded domain with a complete invariant K\"ahler metric
and a finite-volume quotient $\Gamma\quot U$. Then
$\Gamma_0=\Gamma\cap\Aut_0(U)$ is a lattice in $G=\Aut_0(U)$.
\end{lemma}

\begin{proof}
We follow the frame-bundle argument of \cite[p.~94]{FW10}. Let $\mathcal F(U)\to U$ be the unitary frame bundle of the given
metric. Equip $\mathcal F(U)$ with the invariant density obtained
from the volume density on $U$ and Haar probability measure on
the unitary fibers. The group $\Aut(U)$ is a closed subgroup
of the isometry group, and hence a Lie group by the
Myers--Steenrod theorem~\cite{MS39}. Since the metric is complete,
the isometry group acts properly on $U$. Thus the action of
$\Aut(U)$ on $U$, and consequently its lifted action on
$\mathcal F(U)$, is proper. The lifted action is also free:
an isometry fixing a frame fixes its base point and has identity
differential there, so uniqueness of geodesics makes it the
identity on $U$.

Since each fiber has measure one, integration along the fibers gives
\[
    \Vol\bigl(\Gamma\quot\mathcal F(U)\bigr)
       =\Vol(\Gamma\quot U)<\infty.
\]
This equality remains valid when $\Gamma$ has torsion. Indeed,
in a local orbifold chart with finite group $Q$, both quotient
integrals are obtained by dividing the corresponding upstairs
integrals by $|Q|$.

Set $H=G\Gamma$. Since $G$ is normal and open in $\Aut(U)$,
we have $H=\Gamma G$, and $H$ is an open and closed subgroup.
Its action on $\mathcal F(U)$ is free and proper. By the
slice theorem~\cite[Theorem~2.4.1]{DK00}, the orbit projection
\[
             \mathcal F(U)\longrightarrow H\quot\mathcal F(U)
\]
is a locally trivial principal $H$-bundle. Choose a nonempty
coordinate ball $B$ whose closure lies in a trivializing
neighborhood of the orbit space. The inverse image of $B$
can then be identified with $H\times B$, with $H$ acting by
left translation on the first factor. In these coordinates,
invariance of the density gives the product form
\[
                              dh\,d\tau,
\]
where $dh$ is left Haar measure on $H$ and $d\tau$ is a positive
smooth density on $B$. Our choice of $B$ ensures that
$0<\tau(B)<\infty$.

The inclusion $G\hookrightarrow H$ induces a diffeomorphism
\[
    (\Gamma_0\quot G)\times B
       \longrightarrow \Gamma\quot(H\times B),
    \qquad
    (\Gamma_0 g,b)\longmapsto\Gamma(g,b).
\]
Surjectivity follows from $H=\Gamma G$, while injectivity follows
from $\Gamma\cap G=\Gamma_0$. Under this identification, the
quotient density is the product of $d\tau$ and the measure on
$\Gamma_0\quot G$ induced by the left Haar measure $dg=dh|_G$.
Thus the corresponding open subset of
$\Gamma\quot\mathcal F(U)$ has volume
\[
    \Vol_{dg}(\Gamma_0\quot G)\,\tau(B)
       \leq \Vol\bigl(\Gamma\quot\mathcal F(U)\bigr)<\infty.
\]
It follows that $\Gamma_0\quot G$ has finite positive measure.

Let $\bar\mu$ be the induced finite measure on $\Gamma_0\quot G$.
For each $g\in G$, right translation satisfies
\[
                     (R_g)_*\bar\mu=c(g)\bar\mu,
\]
where $c:G\to\mathbb R_{>0}$ is the modular character determined by
this convention. Taking total masses and using
$0<\bar\mu(\Gamma_0\quot G)<\infty$ gives $c(g)=1$ for every $g\in G$.
Thus $G$ is unimodular, and $\Gamma_0$ is a lattice in $G$.
\end{proof}

\subsection{The radical, center and compact factors}

The next two statements exclude compact subgroups normalized by the
quotient group. Their proofs use bounded plurisubharmonic functions,
so the finite-volume Liouville lemma applies directly.

\begin{lemma}\label{lem:finite-normal}
Let $V$ be a contractible bounded domain with a complete invariant
K\"ahler metric and a finite-volume quotient $\Gamma\quot V$.
Every finite subgroup of $\Aut(V)$ normalized by $\Gamma$ is trivial.
\end{lemma}

\begin{proof}
First suppose that $a\in \Aut(V)$ has prime order $p$ and centralizes $\Gamma$.
Since $V$ is a finite-dimensional contractible manifold, Smith
theory~\cite[Chapter~III, Theorem~7.11]{Bre72} implies that the
fixed-point set of $\langle a\rangle$ is nonempty and mod-$p$ acyclic. Choose $x\in V$ with $ax=x$. The function
\[
                u_a(z)=\sup_{\gamma\in\Gamma}
                                  \|\gamma az-\gamma z\|^2
\]
is bounded, locally Lipschitz and plurisubharmonic by
Lemma~\ref{lem:sup}. Since $a$ commutes with $\Gamma$, it is
$\Gamma$-invariant. Lemma~\ref{lem:liouville} makes it constant.
Its value at $x$ is zero, and the term $\gamma=1$ gives $a=1$.

Now let $A<\Aut(V)$ be finite and normalized by $\Gamma$.
Conjugation defines a homomorphism
\[
    \rho:\Gamma\longrightarrow\Aut(A),
    \qquad \rho(\gamma)(b)=\gamma b\gamma^{-1}.
\]
Since $A$ is finite, its automorphism group $\Aut(A)$ is finite.
Thus $\Gamma_1=\ker\rho$ has finite index in $\Gamma$.
By definition, every element of $\Gamma_1$ commutes with every
element of $A$. Moreover, \eqref{eq:index} gives
\[
    \Vol(\Gamma_1\quot V)
       =[\Gamma:\Gamma_1]\,\Vol(\Gamma\quot V)<\infty.
\]
If $A$ were nontrivial, Cauchy's theorem would give an element
$a\in A$ of prime order. This element centralizes $\Gamma_1$,
so the first part of the proof, applied with $\Gamma_1$ in place
of $\Gamma$, would imply $a=1$, a contradiction.
Therefore $A$ is trivial.
\end{proof}

\begin{proposition}\label{prop:compact-normal}
Under the hypotheses of Lemma~\ref{lem:finite-normal}, every compact
subgroup of $\Aut(V)$ normalized by $\Gamma$ is trivial.
\end{proposition}

\begin{proof}
Let $K<\Aut(V)$ be a compact subgroup normalized by $\Gamma$,
and let $K^0$ be its identity component. Its connected center must be
trivial. Indeed, if the torus $Z(K^0)^0$ had positive dimension, then
\[
       Z(K^0)^0[2]\cong(\mathbb Z/2\mathbb Z)^{\dim Z(K^0)^0}
\]
would be a nontrivial finite characteristic subgroup of $Z(K^0)^0$,
and hence would be normalized by $\Gamma$. This contradicts
Lemma~\ref{lem:finite-normal}. Thus $K^0$ is semisimple.

Let
\[
    \mathcal C=\{h\in\Aut(V):hk=kh\text{ for every }k\in K^0\}
\]
be the centralizer of $K^0$. The outer automorphism group of the
semisimple Lie algebra $\mathfrak k$ of $K^0$ is finite
\cite[Chapter~IX, Theorem~5.4]{Hel01}. Since $K^0$ is connected,
differentiation induces an injection
$\operatorname{Out}(K^0)\hookrightarrow\operatorname{Out}(\mathfrak k)$.
Thus $\operatorname{Out}(K^0)$ is
finite, and there is a finite-index subgroup
$\Gamma_1<\Gamma$ whose conjugation action on $K^0$ is by inner
automorphisms. For each $\gamma\in\Gamma_1$,
choose $k\in K^0$ inducing the same conjugation on $K^0$.
Then $k^{-1}\gamma\in \mathcal C$, and hence
\[
                              \Gamma_1\subset K^0\mathcal C.
\]

Suppose that $K^0\ne\{1\}$. Choose an element $a$ of order two
in a maximal torus of $K^0$. Smith
theory~\cite[Chapter~III, Theorem~7.11]{Bre72}, applied to the
finite-dimensional contractible manifold $V$, gives $x\in V$
with $ax=x$. Define
\[
    w(z)=\sup_{h\in \mathcal C}\|h(a(z))-h(z)\|^2,
    \qquad
    W(z)=\max_{k\in K^0}w(kz),
\]
where $\|\cdot\|$ denotes the Euclidean norm.
Equivalently,
\[
    W(z)=\sup_{(h,k)\in \mathcal C\times K^0}
                  \|h(a(kz))-h(kz)\|^2.
\]
Because $V$ is bounded, all the holomorphic differences occurring in
these suprema are uniformly bounded. Lemma~\ref{lem:sup} therefore
shows that both $w$ and $W$ are bounded,
nonnegative, locally Lipschitz and plurisubharmonic; the maximum in
the definition of $W$ exists because $K^0$ is compact and $w$ is
continuous.

Since $a\in K^0$ and every element of $\mathcal C$ commutes with $K^0$,
reindexing the supremum gives $w(cz)=w(z)$ for $c\in \mathcal C$.
Similarly, $W(kz)=W(z)$ for $k\in K^0$, while the commutation of
$\mathcal C$ with $K^0$ gives $W(cz)=W(z)$ for $c\in \mathcal C$. Thus $W$ is invariant
under $K^0\mathcal C$, and hence under $\Gamma_1$. Since
$\Gamma_1\quot V$ has finite volume by \eqref{eq:index},
Lemma~\ref{lem:liouville} makes $W$ constant.

Fix $z\in V$ and choose $k\in K^0$ such that $w(kz)=W(z)$.
Since $w\le W$ and $W$ is constant, $w$ attains its global
maximum at an interior point of $V$. The maximum principle
implies that $w$ is constant. But $ax=x$ gives $w(x)=0$,
so $w=0$. Taking $h=1$ in the definition of $w$ yields
$a(z)=z$ for every $z\in V$, a contradiction. Thus $K^0=\{1\}$.
The compact Lie group $K$ is consequently finite, and
Lemma~\ref{lem:finite-normal} gives $K=\{1\}$.
\end{proof}

\begin{theorem}\label{thm:structure}
If a contractible bounded  $\mathrm{HHR/USq}$ domain $U$ admits a finite-volume
quotient, then $G=\Aut_0(U)$ is centerless, semisimple and has no
compact factors.
\end{theorem}

\begin{proof}
Set $G=\Aut_0(U)$ and $\Gamma_0=\Gamma\cap G$.
The complete K\"ahler--Einstein metric on $U$ is invariant under
$\Aut(U)$, so Lemma~\ref{lem:lattice} shows that $\Gamma_0$
is a lattice in $G$. Moreover, $G$ is normal in $\Aut(U)$,
and consequently $\Gamma_0$ is normal in $\Gamma$. We prove the three assertions in turn.

\smallskip
\underline{Semisimplicity.}
We follow a rank-formula argument. By \cite[Corollary~8.6]{Rag72}, the lattice $\Gamma_0$ has a unique maximal normal solvable subgroup, which we denote by $A$.
For every $\gamma\in\Gamma$, conjugation by $\gamma$ preserves
$\Gamma_0$ and sends $A$ to another maximal normal solvable
subgroup of $\Gamma_0$. Uniqueness therefore gives
\[
                           \gamma A\gamma^{-1}=A.
\]
Thus $A$ is normal in $\Gamma$, not merely in $\Gamma_0$.
Since solvable groups are amenable,
Proposition~\ref{prop:l2-dimension}(i) implies that $A$ is finite.
Lemma~\ref{lem:finite-normal}, applied to this finite subgroup
of $\Aut(U)$ normalized by $\Gamma$, then gives
\[
                                      A=\{1\}.
\]
Let $R=\operatorname{Rad}(G)$ be the solvable radical of $G$,
that is, its largest connected normal solvable subgroup,
and let $K_R$ be a maximal compact subgroup of $R$.
The group $R$ is closed and characteristic in $G$, and
$G/R$ is a connected semisimple Lie group.
Apply Prasad's rank formula~\cite[Lemma~6(2)]{Pra76} to the lattice
$\Gamma_0<G$ and its maximal normal solvable subgroup $A$;
see also the proof of \cite[Proposition~3.3]{FW08}. It gives
\begin{equation}\label{eq:prasad}
    \rank A
       =\dim R-\dim K_R+\rank Z(G/R).
\end{equation}
Here the dimensions are real dimensions, $\rank A$ is the Hirsch rank
occurring in Prasad's formula, and $\rank Z(G/R)$ is the rank of the
discrete abelian group $Z(G/R)$. The subgroup in this formula is
$A$, not $\Gamma_0\cap R$; no assertion that $\Gamma_0\cap R$
is a lattice in $R$ is needed. Since $A=\{1\}$, the left-hand side of
\eqref{eq:prasad} is zero.
Both terms on the right are nonnegative, so
\[
                  \dim R=\dim K_R,
                  \qquad \rank Z(G/R)=0.
\]
The dimension equality implies that $K_R$ is open in $R$.
As $R$ is connected, we obtain $R=K_R$, and hence $R$ is compact. Finally, conjugation by $\Gamma$ preserves $G$ and therefore
its characteristic subgroup $R$. Thus $R$ is a compact subgroup
of $\Aut(U)$ normalized by $\Gamma$.
Proposition~\ref{prop:compact-normal} gives $R=\{1\}$.
This proves that $G$ is semisimple.

\smallskip
\underline{No compact factors.}
Let $\mathfrak g$ be the Lie algebra of $G$, and let
$\mathfrak g_{\mathrm c}$ be the sum of its compact simple ideals.
Denote by $G_{\mathrm c}$ the connected Lie subgroup of $G$
with Lie algebra $\mathfrak g_{\mathrm c}$.
This subgroup is compact: it is the image in $G$ of the
simply connected compact semisimple Lie group with Lie algebra
$\mathfrak g_{\mathrm c}$. Every Lie group automorphism of $G$ permutes the compact
simple ideals of $\mathfrak g$. It therefore preserves
$\mathfrak g_{\mathrm c}$ and its corresponding connected
subgroup $G_{\mathrm c}$.
In particular, conjugation by $\Gamma$ preserves $G_{\mathrm c}$.
Proposition~\ref{prop:compact-normal} now gives
\[
                               G_{\mathrm c}=\{1\}.
\]
Hence $\mathfrak g$ has no compact simple ideals, and $G$
has no compact factors.

\smallskip
\underline{Triviality of the center.}
Since $R=\{1\}$, equation~\eqref{eq:prasad} also gives
\[
                                \rank Z(G)=0.
\]
For a connected semisimple Lie group, $Z(G)$ is discrete. The covering
$G\to G/Z(G)$ gives an exact sequence
\[
             \pi_1(G/Z(G))\longrightarrow Z(G)\longrightarrow 1.
\]
The fundamental group of a connected Lie group is finitely generated,
so $Z(G)$ is a finitely generated abelian group. Its rank is zero, and
it is therefore finite. The center is characteristic in $G$, so it is normalized
by $\Gamma$. Applying Lemma~\ref{lem:finite-normal} once more
yields
\[
                                     Z(G)=\{1\}.
\]
This completes the proof.
\end{proof}

\subsection{The finite-index splitting of
\texorpdfstring{$\Gamma$}{Gamma}}

We now use the structure of $G=\Aut_0(U)$ to obtain a direct-product
decomposition of a finite-index subgroup of $\Gamma$.

\begin{proposition}\label{prop:group-product}
There is a finite-index subgroup $\Gamma_b<\Gamma$ containing
$\Gamma_0=\Gamma\cap G$ such that
\begin{equation}\label{eq:group-product}
    \Gamma_b=\Gamma_0\times\Lambda,
    \qquad
    \Lambda<C_{\Aut(U)}(G),
    \qquad
    \Lambda\cap G=\{1\}.
\end{equation}
Here
\[
    C_{\Aut(U)}(G)
    =
    \{h\in\Aut(U):hg=gh\text{ for every }g\in G\}
\]
is the centralizer of $G$ in $\Aut(U)$, and the product in
\eqref{eq:group-product} is an internal direct product.
\end{proposition}

\begin{proof}
Recall that $G=\Aut_0(U)$ and $\Gamma_0=\Gamma\cap G$.
By Lemma~\ref{lem:lattice} and Theorem~\ref{thm:structure},
$\Gamma_0$ is a lattice in the centerless semisimple group $G$,
which has no compact factors. If $G=\{1\}$, take
$\Gamma_b=\Lambda=\Gamma$.

Assume $G\ne\{1\}$. The outer automorphism group of the semisimple
Lie algebra $\mathfrak g$ of $G$ is finite
\cite[Chapter~IX, Theorem~5.4]{Hel01}. Since $G$ is connected,
differentiation induces an injection
$\operatorname{Out}(G)\hookrightarrow\operatorname{Out}(\mathfrak g)$.
Thus $\operatorname{Out}(G)$ is finite,
and the kernel $\Gamma_a$ of the conjugation homomorphism
$\Gamma\to\operatorname{Out}(G)$ has finite index and contains
$\Gamma_0$. Since $Z(G)=\{1\}$, each $\gamma\in\Gamma_a$
determines a unique $r(\gamma)\in G$ satisfying
\[
    \gamma g\gamma^{-1}
    =r(\gamma)g\,r(\gamma)^{-1}
    \qquad(g\in G).
\]
Uniqueness shows that $r:\Gamma_a\to G$ is a homomorphism
with $r|_{\Gamma_0}=\id$.

Write $N=N_G(\Gamma_0)$. Since $\Gamma_0\triangleleft\Gamma$ and
$r(\gamma)$ induces on $G$ the same conjugation as $\gamma$, we have
$r(\Gamma_a)\subset N$. For each $\delta\in\Gamma_0$, the map
\[
    N^0\longrightarrow\Gamma_0,
    \qquad n\longmapsto n\delta n^{-1},
\]
is continuous and takes values in the discrete group $\Gamma_0$.
It is therefore constant; evaluating at the identity shows that
$N^0\subset C_G(\Gamma_0)$. Since $\Gamma_0$ is a lattice in the
connected semisimple group $G$, which has no compact factors, Borel
density~\cite{Rag72} gives
\[
                    C_G(\Gamma_0)=Z(G)=\{1\}.
\]
Thus $N^0=\{1\}$, and the closed subgroup $N<G$ is discrete. The
natural map
\[
                  \Gamma_0\quot G\longrightarrow N\quot G
\]
is a covering with fiber $N/\Gamma_0$. The quotient $N\quot G$ has
positive Haar measure because $N$ is discrete. Since the covering
preserves the local Haar density and $\Vol(\Gamma_0\quot G)<\infty$,
its fiber is finite. Consequently, $[N:\Gamma_0]<\infty$.

Set
\[
    \Gamma_b=r^{-1}(\Gamma_0),
    \qquad
    \Lambda=\ker r.
\]
Since $\Gamma_0\triangleleft N$, the subgroup $\Gamma_b$ is the
kernel of the homomorphism $\Gamma_a\to N/\Gamma_0$; hence it has
finite index in $\Gamma$. Moreover, since $r|_{\Gamma_0}=\id$,
restriction of $r$ gives the split exact sequence
\[
    1\longrightarrow\Lambda\longrightarrow\Gamma_b
      \xrightarrow{\,r\,}\Gamma_0\longrightarrow1,
    \qquad
    \Gamma_0\hookrightarrow\Gamma_b,
\]
where the section is the inclusion. For $\lambda\in\Lambda$, the
defining identity for $r$ gives $\lambda g\lambda^{-1}=g$ for every
$g\in G$, so $\Lambda<C_{\Aut(U)}(G)$. Furthermore,
\[
    \Lambda\cap\Gamma_0=\{1\},
    \qquad
    \Lambda\cap G
       =\Lambda\cap(\Gamma\cap G)
       =\Lambda\cap\Gamma_0
       =\{1\}.
\]
Thus the two factors commute and have trivial intersection, and the
split extension is the internal direct product
$\Gamma_b=\Gamma_0\times\Lambda$, proving \eqref{eq:group-product}.
\end{proof}

\section{Finite volume of a fixed manifold}\label{sec:fixed-volume}

Assume that $G\ne\{1\}$, and retain
$\Gamma_b=\Gamma_0\times\Lambda$ from
Proposition~\ref{prop:group-product}. We will construct an involution
$a\in G$ whose fixed manifold $E=U^a$ satisfies
\[
                  \Vol(\Lambda\quot E,\ke|_E)<\infty.
\]
The argument uses a uniform tube around $E$ and does not presuppose
a product decomposition of $U$. 

\subsection{Replacing the lattice}

We first show that replacing the lattice in $G$ changes quotient
volume by the ratio of the two covolumes.

\begin{lemma}\label{lem:replacement}
For every lattice $L<G$, the group $L\times\Lambda$ is discrete,
acts properly on $U$, and satisfies
\begin{equation}\label{eq:replacement}
 \Vol((L\times\Lambda)\quot U)
   =\frac{\covol_G(L)}{\covol_G(\Gamma_0)}
                         \Vol((\Gamma_0\times\Lambda)\quot U).
\end{equation}
Both covolumes are computed using the same Haar measure on $G$.
\end{lemma}

\begin{proof}
Since $\Lambda$ centralizes $G$ and $\Lambda\cap G=\{1\}$,
the subgroup $H=G\Lambda$ is topologically $G\times\Lambda$,
where $\Lambda$ has the discrete topology. Indeed, $G$ is open
in the Lie group $\Aut(U)$, and the distinct $\Lambda$-components
lie in distinct cosets of $G$. Thus $H$ is open and closed,
and $L\times\Lambda$ is discrete. Its action is proper because
the $\Aut(U)$-action is proper.

Use the unitary frame bundle $\mathcal F(U)$, with the invariant
measure and unit fiber mass from Lemma~\ref{lem:lattice}. The
$H$-action on this bundle is free and proper. Moreover, $H$ is
unimodular, because $H\cong G\times\Lambda$, $G$ is semisimple, and
$\Lambda$ is discrete. Weil's integration formula for the principal
$H$-bundle
\[
              \mathcal F(U)\longrightarrow H\quot\mathcal F(U)
\]
therefore gives a positive transverse density $d\bar\mu$, independent
of $L$, such that
\[
 \Vol((L\times\Lambda)\quot\mathcal F(U))
   =\Vol((L\times\Lambda)\quot H)
      \int_{H\quot\mathcal F(U)}d\bar\mu.
\]
Since $(L\times\Lambda)\quot H\cong L\quot G$, its volume is
$\covol_G(L)$. Thus, with
\[
 A:=\int_{H\quot\mathcal F(U)}d\bar\mu\in(0,\infty],
\]
we have
\[
 \Vol((L\times\Lambda)\quot\mathcal F(U))=\covol_G(L)A.
\]
Integration over the unitary fibers, which have mass one, identifies
the frame-bundle volumes with the corresponding orbifold volumes on
$U$, including the finite-stabilizer weights. Taking
$L=\Gamma_0$ now shows that $A<\infty$, and taking the ratio proves
\eqref{eq:replacement}.
\end{proof}

The lattices $L$ and $\Gamma_0$ need not be commensurable, so the
replacement quotient need not cover the original one.

\subsection{The Cartan fixed manifold}

We choose a lattice invariant under a Cartan involution and then
use $L^2$ cohomology to realize that involution by an element of $G$.

Choose a Cartan involution $\theta$ of $G$, and let $K=G^\theta$
be its connected maximal compact subgroup. Write
$\frg=\frk\oplus\frp$ for the corresponding decomposition of
Lie algebras. Since $G$ is centerless, its adjoint representation is
faithful; in particular, $G$ is linear. The Borel--Raghunathan
construction, applied to the involution $\sigma=\theta$, gives a
uniform lattice $L'<G$ satisfying $\theta(L')=L'$~\cite{Rag82}; see
also \cite[Theorem~3.1(1)]{Mon23}. Selberg's lemma then supplies a
torsion-free subgroup
$L_1<L'$ of finite index, and
\[
                         L=L_1\cap\theta(L_1)
\]
is torsion-free, uniform and $\theta$-invariant.

Put $\Pi=L\times\Lambda$. Lemma~\ref{lem:replacement} gives
$\Vol(\Pi\quot U)<\infty$. Proposition~\ref{prop:l2-dimension}(i)
and the group $L^2$ K\"unneth formula~\cite{Luc02} yield
\[
 0<\bnum{n}(\Pi)=\sum_{i+j=n}\bnum{i}(L)\bnum{j}(\Lambda)<\infty.
\]
Here the convention is $0\cdot\infty=0$. All summands are
nonnegative, so positivity gives indices $i,j$ for which
$\bnum{i}(L)>0$ and $\bnum{j}(\Lambda)>0$. Since $G/K$ is
contractible and $L$ acts freely and cocompactly, $\bnum{i}(L)$ is
the $i$-th $L^2$-Betti number of $L\quot G/K$. Applying
\cite[Theorem~1.1(a)]{Olb02} therefore gives
\begin{equation}\label{eq:equalrank}
                  \rank_\C\frg_\C=\rank_\C\frk_\C.
\end{equation}

Let $T<K$ be a maximal torus, with Lie algebra $\frt$.
Equation~\eqref{eq:equalrank} makes $\frt_\C$ a Cartan subalgebra
of $\frg_\C$. The complexified involution $\theta_\C$ fixes it
pointwise and acts by scalars on the root spaces. The bracket
relations make these scalars a character of the root lattice.
In the complex adjoint group, this character is evaluation at a
unique element $a$ of the complexified torus $T_\C$. Thus
$\theta_\C=\Int(a)$. Since $\theta_\C^2=1$ and the adjoint
group is centerless, $a^2=1$. Every element of order two in
$T_\C$ belongs to its compact real form $T$. Restriction to $G$
therefore gives
\begin{equation}\label{eq:cartan}
 a\in T,\qquad a^2=1,\qquad \theta=\Int(a),\qquad C_G(a)=K.
\end{equation}
The element $a$ is nontrivial, since $G$ has no compact factors.
It normalizes $L$ and centralizes $\Lambda$.

Set $E=U^a$. Smith theory makes $E$ nonempty and connected,
with vanishing reduced homology over $\mathbf F_2$
\cite[Chapter~III, Theorem~7.11]{Bre72}. Local linearization
of the finite holomorphic action makes $E$ a closed complex
submanifold. It is totally geodesic, hence complete. Its curvature
and all covariant derivatives are bounded by restriction from $U$.
Also $\exp_x^E=\exp_x^U|_{T_xE}$, so
$\operatorname{inj}(E)\ge\operatorname{inj}(U)>0$.
Thus the induced metric has bounded geometry.

The $\Lambda$-action on $E$ is proper and faithful. To check
faithfulness, fix $x\in E$. Its kernel lies in
$\Lambda\cap\operatorname{Stab}_{\Aut(U)}(x)$, which is finite
because the stabilizer is compact and $\Lambda$ is discrete. The
kernel is normal in $\Lambda$ and centralized by $\Gamma_0$,
so it is normalized by $\Gamma_b$. Lemma~\ref{lem:finite-normal}
makes this kernel trivial. It remains to prove finite volume of
$\Lambda\quot E$.

\subsection{Killing fields and displacement}

The next two lemmas turn the adjoint spectral gap of $L$ into a
uniform lower bound for its displacement on $U$.
For a real linear map $A$, let $\rho(A)$ denote the spectral
radius of its complexification.

\begin{lemma}\label{lem:killing-displacement}
Let a positive-dimensional Lie group $G$ act faithfully by
isometries on a connected complete Riemannian manifold $(X,g)$
with $|\mathrm{Rm}|\le\kappa$. There is a constant
$C=C(\dim X,\kappa)>0$ such that
\begin{equation}\label{eq:adjoint-displacement}
 \log\rho(\operatorname{Ad}h)\le C\,d_X(x,hx)
                    \qquad(h\in G,\ x\in X).
\end{equation}
\end{lemma}

\begin{proof}
Identify $\frg$ with its Killing fields and set
\[
                \|Z\|_x^2=|Z(x)|^2+|\nabla Z(x)|^2,
\]
using the Hilbert--Schmidt norm for $\nabla Z$. This is a norm:
a Killing field with zero first jet at $x$ vanishes along every
geodesic from $x$, by the Jacobi equation, and hence vanishes on $X$.

Along a unit-speed geodesic, put $A_Z=\nabla Z$. The Killing
identities give $\nabla_vZ=A_Zv$ and
$\nabla_vA_Z=R(v,Z)$, up to the curvature sign convention.
Hence, for $F=|Z|^2+|A_Z|^2$, the curvature bound gives
$|F'|\le C F$. Equivalently,
$|\frac{d}{dt}\|Z\|_{\gamma(t)}|\le C\|Z\|_{\gamma(t)}$.
Integration along a minimizing geodesic gives
\begin{equation}\label{eq:killing-transport}
                   \|Z\|_y\le e^{C d_X(x,y)}\|Z\|_x.
\end{equation}
Isometries preserve first jets, so
\[
 \|\operatorname{Ad}(h)Z\|_x
   =\|Z\|_{h^{-1}x}\le e^{C d_X(x,hx)}\|Z\|_x.
\]
Bounding the spectral radius by this operator norm proves
\eqref{eq:adjoint-displacement}.
\end{proof}

\begin{lemma}\label{lem:uniform-spectrum}
Let $G\ne\{1\}$ be connected centerless semisimple without
compact factors, and let $L<G$ be a torsion-free uniform lattice.
Then
\begin{equation}\label{eq:lattice-spectrum}
        \eta_L:=\inf_{l\in L\setminus\{1\}}
                       \log\rho(\operatorname{Ad}l)>0.
\end{equation}
For its action on $(U,\ke)$, there is consequently $\delta_L>0$
such that
\begin{equation}\label{eq:lattice-displacement}
 d_{\mathrm{KE}}(x,lx)\ge\delta_L
                    \qquad(x\in U,\ l\in L\setminus\{1\}).
\end{equation}
\end{lemma}

\begin{proof}
Equip $X_0=G/K$ with the metric defined by the Cartan inner
product, and put $o=eK$. Since $L$ is torsion-free, it acts
freely on $X_0$. Compactness of $L\quot X_0$ gives $s_L>0$ with
\[
           d_{X_0}(y,ly)\ge s_L\qquad(y\in X_0,\ l\ne1).
\]
Fix $l\in L$ and let $(y_j)$ minimize its displacement. Choose
$\gamma_j\in L$ so that $z_j=\gamma_jy_j$ lies in a fixed compact
set $Q_0$ with $LQ_0=X_0$, and put
$h_j=\gamma_jl\gamma_j^{-1}$. The numbers
$d_{X_0}(z_j,h_jz_j)$ are bounded. Properness of the $L$-action
therefore leaves only finitely many possibilities for $h_j$.
After passing to a subsequence, $h_j=h$ and $z_j\to z\in Q_0$;
then $h$ attains the minimum displacement of $l$. Conjugating back
shows that $l$ does as well. A nonidentity $l$ is consequently
axial in the Hadamard manifold $X_0$. Its stable translation length
\[
             \tau(l)=\lim_{j\to\infty}j^{-1}d_{X_0}(o,l^jo)
\]
equals this minimum and is at least $s_L$.

If $g=k_1\exp(H)k_2$ is a Cartan decomposition, with $H$ in a
positive chamber of a maximal abelian subspace
$\mathfrak a\subset\frp$, then
\[
 d_{X_0}(o,go)=|H|,\qquad
 \log\|\operatorname{Ad}g\|_{\mathrm{op}}
       =\max_{\alpha\in\Sigma}\alpha(H)\ge c_G|H|.
\]
Here $\Sigma$ is the restricted root system, the operator norm
uses the Cartan inner product, and $c_G>0$ exists because the
roots span $\mathfrak a^*$ and occur in opposite pairs
\cite{Hel01}. Applying this estimate to $l^j$ and using the
spectral-radius formula yields
\[
 \log\rho(\operatorname{Ad}l)
   =\lim_{j\to\infty}j^{-1}\log\|\operatorname{Ad}(l)^j\|_{\mathrm{op}}
   \ge c_G\tau(l)\ge c_Gs_L.
\]
This proves \eqref{eq:lattice-spectrum};
Lemma~\ref{lem:killing-displacement} gives
\eqref{eq:lattice-displacement} with $\delta_L=\eta_L/C$.
\end{proof}

\subsection{A uniform tube and its quotient volume}

We now separate the $L$-translates of $E$ and integrate over a
normal tube that embeds in the finite-volume quotient.

\begin{lemma}\label{lem:fixed-tube}
Let $b$ be an isometric involution of a complete $m$-dimensional
Riemannian manifold $X$ with $|\mathrm{Rm}|\le\kappa$ and
$\operatorname{inj}(X)\ge i_0>0$. Suppose that $E=X^b$ is
nonempty and connected. There is $r_0>0$, depending only on
$m,\kappa,i_0$, such that the normal exponential map identifies
$\{v\in NE:|v|<r\}$ with the metric $r$-neighborhood of $E$
for $0<r<r_0$. Its volume Jacobian is at least $2^{-m}$,
relative to induced volume on $E$ times Euclidean normal-fiber
volume.
\end{lemma}

\begin{proof}
Choose $r_0<i_0/4$ so that $\kappa r_0^2\le1/4$. The fixed
manifold is totally geodesic, and $db_x$ acts as $-1$ on $N_xE$.
If $y=\exp_xv$ with $|v|<r_0$, then
$by=\exp_x(-v)$. The segment from $y$ through $x$ to $by$
has length $2|v|<i_0$, so it is the unique minimizing segment
between its endpoints. Its midpoint recovers $x$, proving
injectivity of the normal exponential map. Conversely, a point
at distance less than $r$ from $E$ has a closest point in $E$,
because Hopf--Rinow makes $X$ proper; a minimizing segment
is normal to $E$. This proves surjectivity onto the metric tube.
More explicitly, on the tube
\[
 x(y)=\exp_y\!\left(\tfrac12\log_y(by)\right),
 \qquad v(y)=\log_{x(y)}y,
\]
and the short logarithms make these maps smooth. They give the
inverse of the normal exponential map.

For the Jacobian, split $T(NE)$ using the normal connection.
A base variation $u\in T_xE$ and a normal variation $w\in N_xE$
give a Jacobi field along $t\mapsto\exp_x(tv)$. Since $E$ is
totally geodesic, its shape operators vanish, and the horizontal
splitting gives the initial data $J(0)=u$, $J'(0)=w$. After
parallel transport to $T_xX$,
the field satisfies
\[
 j(t)=u+tw-\int_0^t(t-s)\mathcal R(s)j(s)\,ds,
                \qquad \|\mathcal R(s)\|\le\kappa|v|^2.
\]
The integral equation gives $\sup_{[0,1]}|j|\le2(|u|+|w|)$ and
\[
 |j(1)-(u+w)|\le\kappa r_0^2(|u|+|w|)
 \le\tfrac12\sqrt{|u|^2+|w|^2}.
\]
Since $u\perp w$, every singular value of the exponential differential
is at least $1/2$, proving the claim.
\end{proof}

\begin{proposition}\label{prop:fixed-volume}
For the Cartan fixed manifold $E=U^a$ constructed above,
\begin{equation}\label{eq:fixed-volume}
                     \Vol(\Lambda\quot E,\ke|_E)<\infty.
\end{equation}
\end{proposition}

\begin{proof}
For $l\in L\setminus\{1\}$, let $b_l=lal^{-1}$. It fixes
$lE$, and
\[
                   b_la=l\theta(l)^{-1}\in L\setminus\{1\}.
\]
The last assertion uses $\theta(L)=L$ and $L\cap K=\{1\}$:
equality $l=\theta(l)$ would put $l$ in $K$.
For $x\in E$ and $y\in lE$, Lemma~\ref{lem:uniform-spectrum}
therefore gives
\[
 \delta_L\le d_{\mathrm{KE}}(x,b_lax)
       =d_{\mathrm{KE}}(x,b_lx)\le2d_{\mathrm{KE}}(x,y).
\]
Consequently,
\begin{equation}\label{eq:translate-separation}
                   d_{\mathrm{KE}}(E,lE)\ge\delta_L/2.
\end{equation}

\begin{figure}[ht]
\centering
\setlength{\unitlength}{1pt}
\begin{picture}(320,82)
  \qbezier(20,18)(145,8)(270,18)
  \qbezier(20,9)(145,1)(270,9)
  \qbezier(20,27)(145,17)(270,27)
  \qbezier(20,62)(145,72)(270,62)
  \qbezier(20,53)(145,63)(270,53)
  \qbezier(20,71)(145,81)(270,71)
  \put(276,14){\makebox(0,0)[l]{$E$}}
  \put(276,58){\makebox(0,0)[l]{$lE$}}
  \put(22,1){\makebox(0,0)[l]{$\mathcal T_r(E)$}}
  \put(22,79){\makebox(0,0)[l]{$\mathcal T_r(lE)$}}
  \put(150,20){\vector(0,1){40}}
  \put(156,40){\makebox(0,0)[l]{$\ge\delta_L/2$}}
\end{picture}
\caption{The radius-$r$ tubes are
disjoint when $2r<\delta_L/2$; no product structure on $U$ is
assumed.}
\label{fig:fixed-tubes}
\end{figure}

Put $m=2n$, $q=m-\dim_\R E$, and choose
$0<r<\min\{r_0,\delta_L/4\}$, where $r_0$ is given by
Lemma~\ref{lem:fixed-tube}. Write $\mathcal T_r(E)$ for the
metric tube. Since
$l\mathcal T_r(E)=\mathcal T_r(lE)$, an intersection between
these two tubes would give $d_{\mathrm{KE}}(E,lE)<2r$, contrary
to \eqref{eq:translate-separation}. Thus the distinct
$L$-translates of $\mathcal T_r(E)$ are disjoint. Since $\Lambda$
preserves the tube, the map
\begin{equation}\label{eq:tube-quotient}
 \Lambda\quot\mathcal T_r(E)\longrightarrow(L\times\Lambda)\quot U
\end{equation}
is an open orbifold embedding: if $l\lambda$ identifies two
points of the tube, disjointness forces $l=1$, so stabilizers
also agree on both sides.

Let $v_q$ be the Euclidean volume of the unit ball in $\R^q$.
Equivariance of the normal exponential map and fiber integration
give
\begin{equation}\label{eq:tube-volume}
 \begin{split}
 \Vol((L\times\Lambda)\quot U)
   &\ge\Vol(\Lambda\quot\mathcal T_r(E))\\
   &\ge2^{-m}v_qr^q\Vol(\Lambda\quot E,\ke|_E).
 \end{split}
\end{equation}
Indeed, on an orbifold chart $V/Q$ of $\Lambda\quot E$, let
$\chi$ be a nonnegative $Q$-invariant cutoff. Normal coordinates,
the fiberwise coarea formula, and Lemma~\ref{lem:fixed-tube}
give the inequality of extended integrals
\[
 \frac1{|Q|}\int_V\chi(x)
       \int_{|v|<r}\!\operatorname{Jac}(\exp^\perp)(x,v)\,dv\,
       d\operatorname{vol}_E(x)
 \ge 2^{-m}v_qr^q\frac1{|Q|}\int_V\chi\,d\operatorname{vol}_E.
\]
Thus the same stabilizer weight occurs on both sides. Summing an
orbifold partition of unity and applying Tonelli proves
\eqref{eq:tube-volume} without assuming that its right-hand side
is finite. Its left-hand side is finite by
Lemma~\ref{lem:replacement}, proving the proposition.
\end{proof}

\begin{rem*}
In the product model $U=S\times Y$, the geodesic symmetry at
$o\in S$ fixes $E=\{o\}\times Y$, and the tube is
$B_S(o,r)\times Y$. Inequality~\eqref{eq:tube-volume} gives the
corresponding volume lower bound before such a product is known.
\end{rem*}

\section{Proof of Theorem~\ref{thm:hhr}}\label{sec:splitting}

Section~\ref{sec:fixed-volume} gives finite quotient volume for the
fixed manifold of a Cartan involution. We first identify this
manifold with the common fixed set of a maximal compact subgroup.
A dimension comparison then turns the orbit map into a smooth
product decomposition. We subsequently identify its complex
structure, metric, and factor actions.

Throughout the intermediate arguments, assume that $n\ge2$ and
$G=\Aut_0(U)\ne\{1\}$. Retain the notation of
Sections~\ref{sec:structure} and~\ref{sec:fixed-volume}. In particular,
$G$ is centerless semisimple without compact factors, and
\[
 \Gamma_b=\Gamma_0\times\Lambda,\qquad
 \Gamma_0=\Gamma\cap G,
\]
where $\Gamma_b$ has finite index in $\Gamma$ and $\Lambda$
centralizes $G$. Let $K<G$ be the connected maximal compact
subgroup, $T<K$ the maximal torus, and $a\in T$ the involution
constructed in Section~\ref{sec:fixed-volume}. For the corresponding
Cartan decomposition,
\[
 \frg=\frk\oplus\frp,\qquad
 \operatorname{Ad}(a)|_{\frp}=-\id,\qquad C_G(a)=K.
\]
We also retain the torsion-free uniform lattice $L<G$ and put
$\Pi=L\times\Lambda$.

The manifold $E=U^a$ is nonempty, connected, and complex. Its
induced K\"ahler metric is complete and has bounded geometry.
The $\Lambda$-action on $E$ is proper and faithful, and
Proposition~\ref{prop:fixed-volume} gives
\[
                    \Vol(\Lambda\quot E,\ke|_E)<\infty.
\]

\subsection{The common fixed manifold}

The torus fixed set supplies both a fixed point and the homology
condition needed for the dimension comparison. We state that
condition explicitly.

\begin{lemma}\label{lem:torus}
Let a compact torus $T$ act smoothly on a nonempty finite-dimensional
manifold $W$ satisfying $\widetilde H_*(W;\Q)=0$. Then
\[
             W^T\ne\varnothing,\qquad
             \widetilde H_*(W^T;\Q)=0.
\]
\end{lemma}

\begin{proof}
For background on torus fixed sets, see
\cite[Chapter~III, \S10]{Bre72}. We give a localization argument
that also applies when $W$ is noncompact.

First suppose that $T=S^1$. Write $H^*_{S^1}$ for Borel
equivariant cohomology, and set
\[
 A=W^{S^1},\qquad B=W\setminus A,\qquad
 R=H^*(BS^1;\Q)=\Q[u],\quad \deg u=2.
\]
The Borel fibration and the homology assumption give
$H^*_{S^1}(W;\Q)\cong R$.

Every stabilizer on $B$ is finite. Consequently, the natural map
\[
             ES^1\times_{S^1}B\longrightarrow B/S^1
\]
induces an isomorphism in rational cohomology. Indeed, over a
slice with finite stabilizer $H$, the local Borel model is
homotopy equivalent to $BH$, whose positive-degree rational
cohomology vanishes by averaging over $H$. The Leray spectral
sequence gives the asserted isomorphism. This uses only local
slice models, not a finite number of orbit types or a uniform
bound on stabilizer orders. The locally free action
of the compact group $S^1$ is proper, so $B/S^1$ is a
finite-dimensional orbifold and has the homotopy type of a CW
complex of dimension at most $\dim W-1$. Its rational cohomology
therefore vanishes in sufficiently high degrees. It follows that
\[
                    H^*_{S^1}(B;\Q)[u^{-1}]=0.
\]
In particular, $A$ cannot be empty, since otherwise this would
contradict $H^*_{S^1}(W;\Q)=R$.

Choose an invariant tubular neighborhood $N$ of $A$. The same
vanishing holds on $N\setminus A$. Localizing the
Mayer--Vietoris sequence for $W=N\cup B$, and using the
equivariant retraction of $N$ onto $A$, yields
\[
 \Q[u,u^{-1}]
   \cong H^*(A;\Q)\otimes_{\Q}\Q[u,u^{-1}].
\]
The module on the right is free over $\Q[u,u^{-1}]$, with rank
$\dim_\Q H^*(A;\Q)$. Thus this dimension is one. Hence $A$ is
connected and its positive-degree rational cohomology vanishes.
The universal coefficient theorem gives
$\widetilde H_*(A;\Q)=0$.

For a general torus, write $T=T_1\times\cdots\times T_r$, with
$T_j\cong S^1$, and set
\[
                  W_0=W,\qquad W_j=W_{j-1}^{T_j}.
\]
The circle case shows inductively that every $W_j$ is nonempty
and rationally acyclic. Since the factors commute, each $T_j$
preserves $W_{j-1}$ and $W_r=W^T$.
\end{proof}

\begin{proposition}\label{prop:common-fixed}
The three fixed sets coincide:
\begin{equation}\label{eq:fixed-equality}
                       Y:=U^K=U^a=U^T.
\end{equation}
Moreover,
\[
                    \widetilde H_*(Y;\Q)=0.
\]
The induced metric on $Y$ is complete and has bounded geometry,
and $\Lambda$ acts properly and faithfully on $Y$ with finite
quotient volume.
\end{proposition}

\begin{proof}
Since $U$ is contractible, Lemma~\ref{lem:torus} supplies a point
$x_0\in U^T\subset E$. Fix $k\in K$. Every element of the
connected compact group $K$ is conjugate into $T$, so there
is $h\in K$ such that $h^{-1}kh\in T$. Thus $k$ fixes $hx_0$.
Moreover, $K=C_G(a)$, so $h$ commutes with $a$ and
\[
                       a(hx_0)=h(ax_0)=hx_0.
\]
Hence $hx_0\in E$. Both $K$ and $\Lambda$ preserve $E$, since
they centralize $a$.

Consider the function
\[
 u_k(z)=\sup_{\lambda\in\Lambda}
                  \|\lambda kz-\lambda z\|^2,\qquad z\in E,
\]
where the norm is Euclidean in $\C^n$. The maps inside the
supremum are holomorphic and uniformly bounded, because $U$
is bounded. Lemma~\ref{lem:sup} therefore makes $u_k$ bounded,
locally Lipschitz, and plurisubharmonic. For $\mu\in\Lambda$,
commutation with $k$ and reindexing give
\[
 u_k(\mu z)
   =\sup_{\lambda\in\Lambda}
        \|\lambda\mu kz-\lambda\mu z\|^2
   =u_k(z).
\]
Hence $u_k$ descends to the complete finite-volume K\"ahler
orbifold $\Lambda\quot E$. It is constant by
Lemma~\ref{lem:liouville}. Its value at $hx_0$ is zero, so
$u_k\equiv0$. Taking $\lambda=1$ shows that $k$ fixes every
point of $E$.

Since $k$ was arbitrary, $E\subset U^K$. The inclusions
\[
                     U^K\subset U^T\subset U^a=E
\]
prove~\eqref{eq:fixed-equality}. Lemma~\ref{lem:torus} now gives
the homology assertion. The metric and action properties are
those already established for $E$.
\end{proof}

\subsection{The smooth product decomposition}

We next show that every $G$-orbit meets $Y$ exactly once.
The dimension comparison supplies local coordinates; a
closedness argument then proves that these coordinates cover $U$.

\begin{lemma}\label{lem:splitting-dimension}
The dimensions satisfy
\begin{equation}\label{eq:dimension}
                  \dim_{\R}(G/K)+\dim_{\R}Y=2n.
\end{equation}
\end{lemma}

\begin{proof}
Equip
\[
                            \mathcal X=(G/K)\times Y
\]
with the product of a $G$-invariant symmetric metric on $G/K$
and the induced metric on $Y$. The manifold $G/K$ is
contractible, so the K\"unneth theorem and
Proposition~\ref{prop:common-fixed} give
$\widetilde H_*(\mathcal X;\Q)=0$. The product metric is complete and
has bounded geometry.

The group $\Pi=L\times\Lambda$ acts by
\[
                   (l,\lambda)(gK,y)=(lgK,\lambda y).
\]
This action is proper and has finite stabilizers. It preserves
orientation: $G$ is connected and $\Lambda$ acts holomorphically
on $Y$. Moreover,
\[
 \Pi\quot \mathcal X
     =(L\quot G/K)\times(\Lambda\quot Y)
\]
has finite volume, since the first factor is compact and the
second has finite volume by
Proposition~\ref{prop:common-fixed}.

Lemma~\ref{lem:replacement} also gives
$\Vol(\Pi\quot U)<\infty$. Thus the same discrete group $\Pi$
acts properly and orientation-preservingly on the two
finite-volume models $U$ and $\mathcal X$, and all hypotheses of
Proposition~\ref{prop:l2-dimension}(ii) are satisfied. It follows
that $\dim_{\R}\mathcal X=2n$, proving \eqref{eq:dimension}.
\end{proof}

\begin{proposition}\label{prop:orbit-diffeomorphism}
The orbit map
\begin{equation}\label{eq:orbit-map}
 \Phi:(G/K)\times Y\longrightarrow U,\qquad
                      \Phi(gK,y)=gy,
\end{equation}
is a diffeomorphism.
\end{proposition}

\begin{proof}
The map is well-defined because $K$ fixes $Y$ pointwise.
For $x\in U$, write
\[
                       G_x=\{g\in G:gx=x\}.
\]
Properness of the action makes $G_x$ compact. If $y\in Y$,
then $K\subset G_y$, and maximal compactness gives $G_y=K$.

We first check injectivity. Recall the Cartan parametrization
\[
              \frp\longrightarrow G/K,\qquad
                         \xi\longmapsto\exp(\xi)K,
\]
which is a diffeomorphism
\cite[Chapter~VI, \S\S1--2]{Hel01}. Since
\[
                       a\exp(\xi)K=\exp(-\xi)K,
\]
the point $o=eK$ is the only fixed point of $a$ on $G/K$.
If $g$ normalizes $K$, then $gK$ is fixed by $K$, and hence
by $a$. Therefore $g\in K$, proving $N_G(K)=K$.
Now suppose that $gy=y'$ for $y,y'\in Y$. Comparing
stabilizers gives
\[
                         gKg^{-1}=K,
\]
so $g\in K$ and $y=y'$. This proves injectivity of $\Phi$.

To compute its differential, let
\[
 \xi_U(x)=\left.\frac{d}{dt}\right|_{t=0}\exp(t\xi)x,
                         \qquad \xi\in\frg,
\]
be the fundamental vector field associated with $\xi$. At
$(o,y)$,
\[
 d\Phi_{(o,y)}(\xi,v)=\xi_U(y)+v,
                    \qquad \xi\in\frp,\quad v\in T_yY.
\]
Evaluation is injective on $\frp$, because its kernel in
$\frg$ is the stabilizer algebra $\frk$. Furthermore,
\[
 da_y\bigl(\xi_U(y)\bigr)
     =\bigl(\operatorname{Ad}(a)\xi\bigr)_U(y)
     =-\xi_U(y),
 \qquad
 da_y(v)=v.
\]
Thus, the subspaces $\frp\cdot y:=\{\xi_U(y):\xi\in\frp\}$ and $T_yY$ meet trivially. Equation~\eqref{eq:dimension} says that their
dimensions add to $2n$, so $d\Phi_{(o,y)}$ is an isomorphism.
Equivariance gives the same conclusion everywhere. Hence
$\Phi$ is a local diffeomorphism and its image is open.

To prove closedness, put $d=\dim_\R(G/K)$. For $x\in U$, choose
$q\in G$ such that $qG_xq^{-1}\subset K$. Then
\[
 \dim Gx=\dim G-\dim G_x\ge \dim G-\dim K=d.
\]
Equality holds if and only if $x\in G\cdot Y$. Indeed, in the
equality case $qG_xq^{-1}$ is a closed subgroup of $K$ with the
same dimension as $K$. It is therefore open in $K$, and hence
equals $K$ because $K$ is connected. Thus $qx\in U^K=Y$.
The converse follows from $G_y=K$ for $y\in Y$.

The orbit dimension is the rank of the smooth bundle map
\[
 \operatorname{ev}:U\times\frg\longrightarrow TU,
                 \qquad (x,\xi)\longmapsto \xi_U(x).
\]
Consequently, the locus on which this rank is at most $d$ is
closed: it is the zero locus of $\bigwedge^{d+1}\operatorname{ev}$.
By the preceding paragraph, this locus is precisely $G\cdot Y$.

The image of $\Phi$ is therefore nonempty, open, and closed
in the connected manifold $U$. It equals $U$, and the
injective local diffeomorphism $\Phi$ is a global
diffeomorphism.
\end{proof}

The two transverse families in the preceding proposition may be
summarized by the diagram
\[
\begin{array}{ccccc}
 (G/K)\times\{y\} & \lhook\joinrel\longrightarrow &
 (G/K)\times Y & \longleftarrow\joinrel\lhook & \{o\}\times Y\\
 {\scriptstyle\sim}\,\downarrow &&
 {\scriptstyle\Phi\;\sim}\,\downarrow &&
 {\scriptstyle\sim}\,\downarrow\\
 Gy & \lhook\joinrel\longrightarrow & U &
 \longleftarrow\joinrel\lhook & Y=U^K .
\end{array}
\]
At a point $y\in Y$, the horizontal and vertical tangent spaces
are respectively the $-1$ and $+1$ eigenspaces of $da_y$.

\subsection{The holomorphic product decomposition}

The involution $a$ separates the two tangent factors by its
eigenvalues. Because $a$ is holomorphic, this separation is
compatible with the complex structure of $U$.

\begin{proposition}\label{prop:holomorphic-splitting}
The symmetric space $G/K$ admits a $G$-invariant complex
structure with bounded symmetric realization $S$ for which
\[
                         \Phi:S\times Y\longrightarrow U
\]
is biholomorphic. Moreover, $Y$ is contractible.
\end{proposition}

\begin{proof}
Let $J$ denote the complex structure of $U$. The preceding
proof gives, for every $y\in Y$,
\[
                    T_yU=(\frp\cdot y)\oplus T_yY.
\]
These are respectively the $-1$ and $+1$ eigenspaces of
$da_y$. Since $da_y$ commutes with $J$, both summands are
complex subspaces. Translating by $G$ shows that every
$G$-orbit has $J$-invariant tangent spaces. The action is proper,
so its orbits are embedded; they are therefore complex
submanifolds of $U$.

For each $y\in Y$, the orbit identification
\[
                    G/K\longrightarrow Gy,\qquad gK\longmapsto gy
\]
induces a $G$-invariant complex structure on $G/K$.
The transformation induced by $a$ fixes $o=eK$, is
holomorphic and isometric, and has differential $-\id$ at
$o$. Its conjugates give the symmetries at all other points.
Thus each of these complex structures makes $G/K$ a
Hermitian symmetric space of noncompact type.

We must check that the complex structure does not depend on $y$.
Decompose the symmetric space into irreducible factors,
\[
                         G/K=X_1\times\cdots\times X_r.
\]
Distinct factors have nontrivial isotropy modules supported on
different factors of $K$, so invariant complex structures preserve
each factor. Thus each $X_i$ is Hermitian and has precisely the two
invariant complex structures $\pm J_i$
\cite[Chapter~VIII, \S\S4--6]{Hel01}. Hence $G/K$ has only
finitely many $G$-invariant complex structures.

Under the evaluation isomorphism
$\operatorname{ev}_y:\frp\to\frp\cdot y$, the structure at $o$ is
\[
        j_y=\operatorname{ev}_y^{-1}\circ
                  J|_{\frp\cdot y}\circ\operatorname{ev}_y.
\]
This endomorphism varies continuously with $y$. Since $Y$ is
connected and the set just described is finite, $j_y$ is
independent of $y$.

Give $G/K$ this common complex structure and identify it
with its bounded symmetric realization $S$
\cite[Chapter~VIII, \S7]{Hel01}. For fixed $y$, the map
$\Phi(\,\cdot\,,y)$ is holomorphic by construction. For
fixed $gK$, the map $\Phi(gK,\,\cdot\,)$ is the restriction
of the holomorphic automorphism $g$ to $Y$. Hence the
differential of the smooth map $\Phi$ is complex-linear
on both tangent factors. It follows that $\Phi$ is
holomorphic. Since its differential is invertible and
$\Phi$ is a diffeomorphism, its inverse is holomorphic as well.

Finally, $S$ is contractible, so the projection
$S\times Y\to Y$ is a homotopy equivalence. The
contractibility of $U$ therefore implies that of $Y$.
\end{proof}

\subsection{Factor automorphisms, metrics, and quotients}

We now identify $U$ with $S\times Y$ through $\Phi$.
The action of $G$ is on $S$ alone. This observation gives
discreteness of $\Aut(Y)$, while the K\"ahler condition
determines the metric decomposition needed later.

\begin{corollary}\label{cor:factor-automorphisms}
The group $\Aut(Y)$ is discrete, and
\[
                              G=\Aut_0(S)
\]
under the product identification.
\end{corollary}

\begin{proof}
Extension by the identity defines a continuous homomorphism
\[
 \iota:\Aut(Y)\longrightarrow\Aut(U),\qquad
                         \iota(\alpha)=\id_S\times\alpha,
\]
for the compact-open topologies. Its image intersects $G$
only in the identity. Indeed, the action of $G$ on $S$ is
effective: an element acting trivially on $S$ acts trivially on
$S\times Y=U$. Elements of $G$ act on $S$ alone, whereas
elements of $\iota(\Aut(Y))$ act on $Y$ alone.
Since $G$ is open in the Lie group $\Aut(U)$,
\[
                         \iota^{-1}(G)=\{\id_Y\}
\]
is open in $\Aut(Y)$. Thus $\Aut(Y)$ is discrete.
Restriction identifies $G$ with a connected subgroup of
$\Aut(S)$, so $G\subset\Aut_0(S)$. Conversely, extension by
$\id_Y$ sends the connected group $\Aut_0(S)$ into
$\Aut_0(U)=G$. These inclusions prove the asserted equality.
\end{proof}

\begin{proposition}\label{prop:factor-metric-volume}
The complete K\"ahler--Einstein metric splits as
\begin{equation}\label{eq:metric-product}
 \ke=g_S\oplus g_Y,\qquad
             \Ric(g_S)=-g_S,\qquad \Ric(g_Y)=-g_Y,
\end{equation}
where both factor metrics are complete. Moreover, as complex
orbifolds,
\[
 \Gamma_b\quot U
       \cong(\Gamma_0\quot S)\times(\Lambda\quot Y),
\]
and
\[
 \Vol(\Gamma_b\quot U,\ke)
   =\Vol(\Gamma_0\quot S,g_S)\,
      \Vol(\Lambda\quot Y,g_Y).
\]
Both factors on the right have finite volume.
\end{proposition}

\begin{proof}
At $(o,y)$, the isometry $a$ acts by $-1$ on the tangent
space to $S$ and by $+1$ on the tangent space to $Y$.
These spaces are therefore orthogonal. The $G$-action is
transitive on $S$ and trivial on $Y$, so the mixed metric
terms vanish everywhere and the vertical metric is
independent of the point of $S$.

Let $\omega$ be the K\"ahler form of $\ke$. We can consequently
write
\[
                           \omega=\omega_S(y)+\omega_Y,
\]
where $\omega_Y$ is pulled back from $Y$ and $\omega_S(y)$
is a family of horizontal two-forms. Let $V_1,V_2$ be product
lifts of local vector fields on $S$, and let $W$ be the product
lift of a local vector field on $Y$. Then $[V_i,W]=0$. Since
all mixed components of $\omega$ vanish and $[V_1,V_2]$ is
horizontal, the identity $d\omega=0$ gives
\[
        0=d\omega(V_1,V_2,W)
          =W\bigl(\omega_S(y)(V_1,V_2)\bigr).
\]
Thus $\omega_S(y)$ is independent of $y$, and
\[
 \omega=\operatorname{pr}_S^*\omega_S
          +\operatorname{pr}_Y^*\omega_Y.
\]
Equivalently, the metric is a product.

In product holomorphic coordinates, the metric matrix is block diagonal,
with determinant $\det(g_S)\det(g_Y)$. Hence the Ricci forms satisfy
\[
                         \rho_{\ke}=\rho_{g_S}+\rho_{g_Y}.
\]
Comparing the horizontal and vertical components of
$\rho_{\ke}=-\omega$ proves \eqref{eq:metric-product}. Finally,
if a maximal geodesic of either factor had finite existence time,
keeping the other coordinate fixed would produce an incomplete
geodesic of $(U,\ke)$. Both factor metrics are therefore complete.

For the quotient statement, $\Gamma_0<G$ acts on $S$ alone.
An element $\lambda\in\Lambda$ preserves $Y=U^K$ because it
centralizes $K$. Since it also commutes with $G$,
\[
           \lambda\Phi(gK,y)
                 =g\lambda y
                 =\Phi(gK,\lambda y).
\]
Thus $\Lambda$ acts on $Y$ alone. Under $\Phi$, the action is
therefore
\[
             (\gamma,\lambda)(s,y)=(\gamma s,\lambda y).
\]
The internal direct product $\Gamma_b=\Gamma_0\times\Lambda$
gives the asserted orbifold product. The volume formula follows
from product integration. In local quotient charts, the
stabilizer is the product of
the two factor stabilizers, so its reciprocal order is
the product of their reciprocal orders. The same formula
therefore holds with orbifold volume.

Finally, $g_S$ is $G$-invariant, and its volume form is a
positive multiple of quotient Haar measure on $G/K$.
Since $\Gamma_0$ is a lattice in $G$ and $K$ is compact,
$\Gamma_0\quot S$ has positive finite volume. The metric
$g_Y$ is precisely the induced metric on $Y$, whose quotient
volume is finite by Proposition~\ref{prop:common-fixed}.
A point factor has volume one; if $Y$ is a point,
faithfulness of its $\Lambda$-action gives $\Lambda=\{1\}$.
\end{proof}

\begin{proof}[Proof of Theorem~\ref{thm:hhr}]
If $n=1$, uniformization gives $U\cong\D$, and we take
$S=\D$ and $Y$ to be a point. Suppose $n\ge2$. If
$G=\{1\}$, then $\Aut(U)$ is discrete, and the assertion
holds with $S$ a point and $Y=U$. If $G\ne\{1\}$, Proposition~\ref{prop:holomorphic-splitting}
gives a biholomorphism $U\cong S\times Y$, where $S$ is a
bounded symmetric domain and $Y$ is contractible.
Corollary~\ref{cor:factor-automorphisms} shows that
$\Aut(Y)$ is discrete. This proves the theorem.
\end{proof}

\begin{rem*}
In the proof of Theorem~\ref{thm:hhr}, HHR/USq is used to obtain
the complete K\"ahler--Einstein metric with bounded geometry.
Thus, the same conclusion holds for a contractible bounded
pseudoconvex domain whose complete K\"ahler--Einstein metric has
bounded geometry and admits a finite-volume quotient. Here, bounded
geometry includes a positive lower bound on the injectivity radius,
not merely completeness. 
\end{rem*}

\section{Convex boundary points}\label{sec:convex-geometry}

We now turn to Theorem~\ref{thm:convex}. The argument requires a boundary point that
is complex extreme and whose normal cone has totally real span.
Complex extremality supplies a global holomorphic peak function;
finite quotient volume then makes the point an orbit accumulation
point. Total reality supplies the tangential direction
used to construct a nonzero complete holomorphic vector field in
Section~\ref{sec:flows}. We first find a point with both properties by a purely convex
argument; the quotient enters only afterwards, in Lemma~\ref{lem:peak}.

In this section, $\mathcal K$ denotes a compact convex body in
$\C^n=\R^{2n}$ with nonempty interior, and $n\ge1$. We use the
real Euclidean inner product and write $Jv=iv$. Define the support
function, normal cone and its real span by
\[
 \begin{split}
  h_{\mathcal K}(\eta)&=\max_{x\in \mathcal K}\ip{\eta}{x},\\
  N_{\mathcal K}(p)&=\{\eta:\ip{\eta}{x-p}\le0\text{ for every }x\in \mathcal K\},\\
  V_p&=\Span_\R N_{\mathcal K}(p).
\end{split}
\]
Here $N_{\mathcal K}(p)$ is defined for every $p\in\mathcal K$; it is
$\{0\}$ when $p$ is an interior point.  A real subspace $V$ is
\emph{totally real} if $V\cap JV=\{0\}$.
A point $p\in\partial \mathcal K$ is \emph{complex extreme} if $\mathcal K$ contains
no nonconstant complex affine disk centered at $p$. Let $\CE(\mathcal K)$
denote the set of these points. For later use, the complex-linear
functional $\ell_\eta(z)=\sum z_\alpha\overline{\eta_\alpha}$ satisfies
\begin{equation}\label{eq:linear}
                       \ell_\eta(v)=\ip{\eta}{v}+i\ip{J\eta}{v}.
\end{equation}

\subsection{An integral over the normal bundle}

A positive normal-bundle integral will locate a boundary point
with the two geometric properties required for the flow argument.

\begin{proposition}\label{prop:normal-point}
Every compact convex body $\mathcal K\subset\C^n$ with nonempty interior
has a point $p\in\CE(\mathcal K)$ for which $V_p$ is totally real.
\end{proposition}

Let $P_{\mathcal K}$ be nearest-point projection onto $\mathcal K$. On $\C^n\setminus \mathcal K$
define
\[
                        \mathbf n(y)
                            =\frac{y-P_{\mathcal K}(y)}{\dist(y,\mathcal K)}.
\]
For $r>0$, the parallel surface and its normal parametrization are
\[
        \Sigma_r=\partial(\mathcal K+r\overline{\B^n}),\qquad
        \Psi_r(y)=\left(P_{\mathcal K}(y),\frac{y-P_{\mathcal K}(y)}r\right).
\]
The surface $\Sigma_r$ is $C^{1,1}$. The map $\Psi_r$ is a bi-Lipschitz
homeomorphism onto the unit normal bundle
$\{(p,\nu):p\in\partial \mathcal K,\ \nu\in N_{\mathcal K}(p),\ |\nu|=1\}$,
with inverse $(p,\nu)\mapsto p+r\nu$. On the ambient product space
put
\[
           \alpha=\ip{J\nu}{dp},\qquad
           \omega=\alpha\wedge(d\alpha)^{n-1},\qquad
           I_r(\mathcal K)=\int_{\Sigma_r}\Psi_r^*\omega.
\]
The surface has outward orientation, and the pullback uses
almost-everywhere derivatives. The role of this integral is to find
a normal pair $(p,\nu)$ at which the pulled-back top-degree form is
nonzero. We will prove that $I_r(K)>0$ and choose a parallel surface
on which the nearest-point projection is differentiable almost
everywhere. At a point of nonvanishing, either of the two unwanted
complex directions would force a contradiction: a nonzero complex
subspace of $V_p$, or a nonconstant complex affine disk in $K$
centered at $p$, would produce a nonzero tangent vector annihilating
both $\Psi_r^*\alpha$ and $\Psi_r^*d\alpha$.
Their top-degree wedge would then vanish. The next lemma establishes
the positivity and continuity of the integral, and the proof of
Proposition~\ref{prop:normal-point} carries out these two pointwise arguments.

\begin{lemma}\label{lem:normal-integral}
For fixed $r>0$, $I_r$ is continuous under Hausdorff convergence of
convex bodies with uniform positive inner radius and bounded outer
radius. Moreover, $I_r(\mathcal K)>0$ for every full-dimensional convex body.
\end{lemma}

\begin{proof}
Translate so that $0\in\operatorname{int}\mathcal K$; this does not change
the integral. First suppose $\mathcal K$ is smooth and strongly convex. Its
inverse Gauss map is $p(\nu)=\nabla h_{\mathcal K}(\nu)$.  Let
$\iota(\nu)=(p(\nu),\nu)$ and put
$q(\nu)=\ip{J\nu}{p(\nu)}$.  All forms in the next computation are
pulled back by $\iota$. With $d^cu=-du\circ J$, we have
\begin{equation}\label{eq:normal-dc}
                     \iota^*\alpha=d^ch_{\mathcal K}+dq,
             \qquad \iota^*d\alpha=dd^ch_{\mathcal K}.
\end{equation}
Indeed,
$dq=\ip{J\,d\nu}{p}+\ip{J\nu}{dp}$, whereas
$d^ch_{\mathcal K}=-\ip{J\,d\nu}{p}$.
The map $\nu\mapsto\nabla h_{\mathcal K}(\nu)+r\nu$ preserves outward
orientation. The exact term in \eqref{eq:normal-dc} integrates to
zero after wedging with $(dd^ch_{\mathcal K})^{n-1}$. Hence
\begin{equation}\label{eq:normal-ma}
 I_r(\mathcal K)=\int_{\partial\B^n}d^ch_{\mathcal K}\wedge(dd^ch_{\mathcal K})^{n-1}
                              =\int_{\B^n}(dd^ch_{\mathcal K})^n.
\end{equation}
For the last equality, apply Stokes' theorem on
$\B^n\setminus t\overline{\B^n}$ and let $t\downarrow0$.
Homogeneity makes the inner boundary integral $t^n$ times the outer
one.  Equivalently, the Bedford--Taylor measure~\cite{BT82} of the locally bounded
plurisubharmonic function $h_{\mathcal K}$ does not charge the singleton
$\{0\}$. This proves \eqref{eq:normal-ma}.

We next prove continuity without a smoothness assumption. Suppose
$\mathcal K_j\to \mathcal K$. After translating a common interior
center to the origin and discarding finitely many terms, fix
$0<a<R$ such that
\[
 \B_a^n\subset\mathcal K_j\subset\B_R^n,
 \qquad
 \B_a^n\subset\mathcal K\subset\B_R^n .
\]
Write the boundary of
$\mathcal K_j+r\overline{\B^n}$ as
$\varphi_j(\zeta)=\varrho_j(\zeta)\zeta$, $\zeta\in\partial\B^n$.
The common radius bounds and Hausdorff convergence give
\[
 \varrho_j\longrightarrow\varrho\quad\text{uniformly},
 \qquad \sup_j\operatorname{Lip}(\varphi_j)<\infty.
\]
Nearest-point projections converge uniformly on bounded sets: every
subsequential limit of projected points lies in $\mathcal K$ and is a
nearest point there, and uniqueness identifies the limit.  Hence
\[
 \Xi_j(\zeta)=\left(P_{\mathcal K_j}(\varphi_j(\zeta)),
       \frac{\varphi_j(\zeta)-P_{\mathcal K_j}(\varphi_j(\zeta))}{r}\right)
\]
converges uniformly to the analogous map $\Xi$, while the
$1$-Lipschitz property of nearest-point projection gives
$\sup_j\operatorname{Lip}(\Xi_j)<\infty$.

Because the origin lies in the interior, every radial map $\varphi_j$ is
orientation preserving. Consequently,
\[
        \int_{\partial\B^n}\Xi_j^*\omega=I_r(\mathcal K_j),
        \qquad
        \int_{\partial\B^n}\Xi^*\omega=I_r(\mathcal K).
\]

Apply Lipschitz Stokes to the straight homotopy
$G_j(t,\nu)=(1-t)\Xi(\nu)+t\Xi_j(\nu)$ in the ambient product.
It gives
\[
           I_r(\mathcal K_j)-I_r(\mathcal K)
                 =\int_{[0,1]\times\partial\B^n}G_j^*d\omega.
\]
At almost every point, every top-degree term contains the factor
$\partial_tG_j=\Xi_j-\Xi$; the remaining $2n-1$ factors are bounded
by the common Lipschitz constant.  The images lie in one compact
subset of the ambient product; therefore
\begin{equation}\label{eq:normal-continuity}
                   |I_r(\mathcal K_j)-I_r(\mathcal K)|
                        \le C\|\Xi_j-\Xi\|_\infty\longrightarrow0.
\end{equation}
This is the standard Stokes formula for Lipschitz maps, equivalently
for the associated integral currents; the form $\omega$ is smooth on
the ambient product.

It remains to obtain a positive lower bound that survives
approximation. Choose smooth strongly convex approximations
\cite{Sch14} with $\B_a^n\subset \mathcal K_j\subset \B_R^n$, where
$0<a<R$. Fix $0<\varepsilon<c<a$ and put
\[
 q(z)=c+\varepsilon(|z|^2-1),\qquad
 \rho=\frac{c-\varepsilon}{2R},\qquad
 A_j=\{z\in\B^n:h_{\mathcal K_j}(z)<q(z)\}.
\]
If $|z|<\rho$, then
$h_{\mathcal K_j}(z)\le R|z|<(c-\varepsilon)/2<q(z)$, so
$\B_\rho^n\subset A_j$. On $\partial\B^n$ we have
$h_{\mathcal K_j}\ge a>c=q$. Thus the continuous plurisubharmonic
functions $h_{\mathcal K_j}$ and $q$ satisfy the boundary inequality
required by the Bedford--Taylor comparison principle, and
$A_j\Subset\B^n$. By that principle~\cite{BT82} and
\eqref{eq:normal-ma},
\begin{equation}\label{eq:positive-mass}
 I_r(\mathcal K_j)\ge\int_{A_j}(dd^ch_{\mathcal K_j})^n
       \ge\int_{A_j}(dd^cq)^n
       \ge\int_{\B_\rho^n}(dd^cq)^n>0.
\end{equation}
The lower bound is independent of $j$, so continuity gives
$I_r(\mathcal K)>0$.
\end{proof}

\begin{proof}[Proof of Proposition~\ref{prop:normal-point}]
The projection $P_{\mathcal K}$ is the gradient of the convex $C^{1,1}$ function
\[
                   Q(y)=\sup_{x\in \mathcal K}
                                 \bigl(\ip{y}{x}-\tfrac12|x|^2\bigr).
\]
Its ambient derivative exists and is symmetric almost everywhere.
Let $\mathcal E\subset\C^n\setminus\mathcal K$ be the null set on which this
derivative does not exist as a symmetric map. Since
$|\nabla\dist(\,\cdot\,,\mathcal K)|=1$ almost everywhere outside
$\mathcal K$, the coarea formula gives
\[
 0=\int_{\mathcal E}|\nabla\dist(y,\mathcal K)|\,dy
   =\int_0^\infty\mathcal H^{2n-1}(\mathcal E\cap \Sigma_r)\,dr.
\]
Choose $r>0$ for which $\mathcal H^{2n-1}(\mathcal E\cap \Sigma_r)=0$.
Since $I_r(\mathcal K)>0$, there is a point
$y\in \Sigma_r\setminus \mathcal E$ at which $\Psi_r^*\omega\ne0$.
Write
\[
       p=P_{\mathcal K}(y),\qquad \nu=(y-p)/r,\qquad
       A=DP_{\mathcal K}(y),\qquad
       B=D\mathbf n(y)=\frac{I-A-\nu\otimes\nu}{r}.
\]
The matrices $A$ and $B$ are symmetric, and $A\nu=0$ because
projection is constant along the normal ray. On
$T_y\Sigma_r=\nu^\perp$, put $\beta=\Psi_r^*\alpha$ and
$\tau=\Psi_r^*d\alpha$. Then
\begin{equation}\label{eq:normal-forms}
 \begin{split}
  \beta(u)&=\ip{J\nu}{Au},\\
  \tau(u,v)&=\ip{JBu}{Av}-\ip{JBv}{Au},
                     \qquad \beta\wedge\tau^{n-1}\ne0.
 \end{split}
\end{equation}

We first show that $V_p$ is totally real. If $\eta\in N_{\mathcal K}(p)$
and $t\ge0$, then $r\nu+t\eta\in N_{\mathcal K}(p)$, and hence
$P_{\mathcal K}(y+t\eta)=p$. Taking the one-sided derivative gives $A\eta=0$.
Therefore $V_p\subset\ker A$. If $W=V_p\cap JV_p\ne\{0\}$, then
$W$ is $J$-invariant and has real dimension at least two. Hence
$W\cap\nu^\perp$ contains a nonzero vector $w$.
We have $Aw=AJw=0$ and $Bw=w/r$. By symmetry of $A$,
\[
                    \beta(w)=0,\qquad
                    \tau(w,u)=r^{-1}\ip{Jw}{Au}=0
\]
for every tangent vector $u$. Thus
$\iota_w(\beta\wedge\tau^{n-1})=0$, which is impossible for a
nonzero top form and a nonzero tangent vector.
This contradiction proves total reality.

Suppose now that a complex disk centered at $p$ lies in $\mathcal K$, and
let $\mathcal L$ be its complex direction line. In particular, $J\mathcal L=\mathcal L$ and
$\mathcal L\subset\nu^\perp$. The supporting hyperplane with
normal $\nu$ contains the disk. Consequently $\nu$ is normal to
$\mathcal K$ at both $p+tv$ and $p+tJv$ for all sufficiently small
real $t$. Thus
\[
 P_{\mathcal K}(y+tv)=p+tv,\qquad
 P_{\mathcal K}(y+tJv)=p+tJv.
\]
It follows that $Av=v$, $AJv=Jv$, and $Bv=BJv=0$. For nonzero $v\in \mathcal L$,
formula~\eqref{eq:normal-forms} gives
\[
                  \beta(v)=0,\qquad
                  \tau(v,u)=-\ip{JBu}{v}=\ip{u}{BJv}=0.
\]
This is the same contradiction. Hence $p\in\CE(\mathcal K)$.
When $n=1$, the top form is simply $\beta$, and its nonvanishing
gives the conclusion in the same way. The choice of $r$ by coarea
ensures that the ambient derivatives used above exist on the chosen
parallel surface.
\end{proof}

\subsection{Peak points and orbit limits}
A global peak function allows us to detect convergence to a
specified boundary point. A boundary point $p$ of a bounded domain $U$ is a \emph{global holomorphic
peak point} if there is $f\in\mathcal O(U)\cap C(\overline U)$ such that
\[
 f(p)=1,\qquad |f(z)|<1\quad(z\in\overline U\setminus\{p\}).
\]
We call $f$ a \emph{peak function at $p$}. For every neighborhood $V$ of
$p$, compactness gives $\sup_{\overline U\setminus V}|f|<1$.
Consequently, $|f(z_j)|\to1$ forces $z_j\to p$.
For $o\in U$, an \emph{orbit accumulation point} is a limit of a sequence
$\gamma_j o$, with $\gamma_j\in\Gamma$.

\begin{lemma}\label{lem:peak}
Let a bounded domain $U$ carry a complete invariant K\"ahler metric
with a finite-volume quotient $\Gamma\quot U$. Every peak point for
$\mathcal O(U)\cap C(\overline U)$ is a limit of $\Gamma o$ for
every $o\in U$.
\end{lemma}

\begin{proof}
Let $p$ be the peak point, and choose $f$ with $f(p)=1$ and
$|f|<1$ on $\overline U\setminus\{p\}$. The function
\[
                         u(z)=\sup_{\gamma\in\Gamma}|f(\gamma z)|^2
\]
is constant by Lemmas~\ref{lem:liouville} and~\ref{lem:sup}.
Since $|f|^2\le u\le1$ and $\sup_U|f|=1$, its value is one.
Choose $\gamma_j$ with $|f(\gamma_jo)|\to1$. Compactness of
$\overline U$ and uniqueness of the peak imply $\gamma_jo\to p$.
\end{proof}

For a bounded convex domain, every complex extreme point is a global
peak point~\cite[Theorem~6]{NP03}. Indeed, an affine complex disk in
$\overline U$ centered at a boundary point lies in every supporting
hyperplane through its center. Thus complex extremality is precisely
the vanishing of the complex lineality set denoted by $L(p)$ in that
theorem. The point supplied by
Proposition~\ref{prop:normal-point} belongs to every orbit limit set.
The use of the normal bundle is relevant even for the bidisk: its
complex extreme set is $\partial\D\times\partial\D$, which has zero
hypersurface measure. At $(1,1)$, the normal cone is
$\R_{\ge0}e_1+\R_{\ge0}e_2$ and its span is totally real. The
argument requires no positive boundary measure for the complex
extreme set.

\section{Complete holomorphic flows and Theorem~\ref{thm:convex}}\label{sec:flows}

Let $U\subset\C^n$ be bounded and convex, and put $\mathcal K=\overline U$.
We will rescale translations in a direction tangent to all supporting
hyperplanes at a suitable boundary orbit limit. The main point is to
prove completeness of the limiting field. For that purpose we need
estimates that hold for all supporting directions at once.

\subsection{Supporting gaps and the Kobayashi metric}

In this section $d_U^{\mathrm K}$ is the Kobayashi distance, and $k_U$ is the
infinitesimal Kobayashi metric normalized by $k_\D(0;v)=|v|$.
Write
\[
           \delta(z)=\dist(z,\partial U),\qquad
           s_\theta(z)=h_{\mathcal K}(\theta)-\ip{\theta}{z},\quad |\theta|=1.
\]
Thus $s_\theta$ is the Euclidean gap from a supporting hyperplane.

\begin{lemma}\label{lem:gaps}
For $x,y,z\in U$ and $|\theta|=1$,
\begin{align}
 e^{-2d_U^{\mathrm K}(x,y)}s_\theta(x)&\le s_\theta(y)
                       \le e^{2d_U^{\mathrm K}(x,y)}s_\theta(x),\label{eq:gap-harnack}\\
 \delta(z)&=\min_{|\theta|=1}s_\theta(z),\label{eq:depth}\\
 \delta(y)&\ge e^{-2d_U^{\mathrm K}(x,y)}\delta(x).\label{eq:depth-harnack}
\end{align}
Kobayashi balls are relatively compact. If $U\subset \B_R^n$, then
\begin{equation}\label{eq:kob-bounds}
                      \frac{|\xi|}{R}\le k_U(z;\xi)
                                               \le\frac{|\xi|}{\delta(z)}.
\end{equation}
\end{lemma}

\begin{proof}
The holomorphic map $h_{\mathcal K}(\theta)-\ell_\theta$ takes $U$ into the
right half-plane $\mathbb H=\{\operatorname{Re}w>0\}$. In our
normalization its metric is $|dw|/(2\operatorname{Re}w)$, so
\[
              |d\log\operatorname{Re}w|
                   \le\frac{|dw|}{\operatorname{Re}w}
                   =2k_{\mathbb H}(w;dw).
\]
Integrate along paths, take the infimum and apply Schwarz--Pick.
This gives
$|\log s_\theta(y)-\log s_\theta(x)|\le2d_U^{\mathrm K}(x,y)$, proving
\eqref{eq:gap-harnack}.

The body $\mathcal K$ is the intersection of its supporting half-spaces.
The largest ball centered at $z$ and contained in $\mathcal K$ therefore
has radius $\min s_\theta(z)$, which proves \eqref{eq:depth}.

Minimizing \eqref{eq:gap-harnack} gives \eqref{eq:depth-harnack}.
A Kobayashi ball of radius $T$ about $x$ consequently stays at
Euclidean distance at least $e^{-2T}\delta(x)$ from the boundary.
Since $U$ is bounded, the ball is relatively compact.

The inclusion $B(z,\delta(z))\subset U$ gives the upper bound in
\eqref{eq:kob-bounds}. For the lower bound, apply Schwarz--Pick to
$\ell_e/R:U\to\D$, where $e=\xi/|\xi|$ if $\xi\ne0$.
\end{proof}

\subsection{Rescaling a tangential translation}

We separate the boundary estimate, the construction of a nonzero
limiting field, and the completeness argument.

\begin{proposition}\label{prop:flow}
Suppose that $\gamma_j\in\Aut(U)$ and
$\gamma_j(o)\to p\in\partial U$, and that $V_p$ is totally real.
Then $U$ has a nowhere-zero holomorphic vector field whose real
flow is complete.
\end{proposition}

Throughout this subsection, assume the hypotheses of
Proposition~\ref{prop:flow}. Since $U$ is bounded, Montel's theorem
allows us to pass to a subsequence such that
\[
\gamma_j\longrightarrow f
\quad\text{locally uniformly on }U,
\qquad
f(o)=p,
\qquad
f(U)\subset\mathcal K.
\]
Here $f:U\to\C^n$ is holomorphic and
$\mathcal K=\overline U$. Put
\[
p_j=\gamma_j(o),
\qquad
\delta_j=\delta(p_j)\longrightarrow0,
\qquad
V=V_p.
\]

\begin{lemma}[Tangential translation estimates]
\label{lem:flow-segments}
Fix $v\in V^\perp$, $z_*\in U$, and $T>0$, and put
\[
 x_j=\gamma_j(z_*),\qquad
 c=e^{-2d_U^{\mathrm K}(o,z_*)}.
\]
For all sufficiently large $j$,
\begin{equation}\label{eq:segment}
s_\theta(x_j+t\delta_jv)\ge\frac{c\delta_j}{2}
\qquad (|t|\le T,\ |\theta|=1).
\end{equation}
Consequently, the whole translated segment lies in $U$ and
has Euclidean distance at least $c\delta_j/2$ from its boundary $\partial U$.
\end{lemma}

\begin{proof}
We first identify the normal cones along the limiting image.
For each fixed $\eta$, the function
\[
z\longmapsto
h_{\mathcal K}(\eta)-\operatorname{Re}\ell_\eta(f(z))
\]
is nonnegative and pluriharmonic. By the minimum principle, it
vanishes somewhere if and only if it vanishes identically. For
$q\in\mathcal K$, the equality
$h_{\mathcal K}(\eta)=\ip{\eta}{q}$ is equivalent to
$\eta\in N_{\mathcal K}(q)$. Evaluating the pluriharmonic function
at $z$ and at $o$ therefore gives
\begin{equation}\label{eq:equal-cones}
N_{\mathcal K}(f(z))=N_{\mathcal K}(p),
\qquad z\in U.
\end{equation}

By \eqref{eq:gap-harnack} and \eqref{eq:depth},
\[
s_\theta(x_j)\ge c\,s_\theta(p_j)\ge c\delta_j
\qquad (|\theta|=1).
\]
Suppose that \eqref{eq:segment} fails along a subsequence.
Choose $t_j\in[-T,T]$ and unit vectors $\theta_j$ witnessing
the failure. The affine identity
\[
s_{\theta_j}(x_j+t_j\delta_jv)
=s_{\theta_j}(x_j)
-t_j\delta_j\ip{\theta_j}{v}
\]
then gives
\[
t_j\ip{\theta_j}{v}>\frac c2,
\qquad
0<s_{\theta_j}(x_j)
<
\left(\frac c2+T|v|\right)\delta_j.
\]
After passing to a subsequence, $\theta_j\to\theta$.
Since $x_j\to f(z_*)$ and $\delta_j\to0$, continuity of the
support function implies $s_\theta(f(z_*))=0$. Thus
\[
\theta\in N_{\mathcal K}(f(z_*))
=N_{\mathcal K}(p)\subset V.
\]
It follows that $\ip{\theta_j}{v}\to0$, contradicting the
first inequality because $|t_j|\le T$. Finally, \eqref{eq:segment} and the supporting half-space
description of $\mathcal K$ place a Euclidean ball of radius
$c\delta_j/2$ about each translated point inside $\mathcal K$.
This proves both membership in $U=\operatorname{int}\mathcal K$
and the asserted boundary-distance estimate.
\end{proof}

\begin{lemma}[A nonzero limit of rescaled fields]
\label{lem:flow-limit}
There exist a nonzero vector $v\in V^\perp$ and a subsequence
such that the holomorphic vector fields
\begin{equation}\label{eq:fields}
Z_j(z)=\delta_j[D\gamma_j(z)]^{-1}v
\end{equation}
converge locally uniformly on $U$ to a nowhere-zero
holomorphic vector field $Z$.
\end{lemma}

\begin{proof}
Choose unit vectors $\nu_j$ with
$s_{\nu_j}(p_j)=\delta_j$. After passing to a subsequence,
$\nu_j\to\nu$ with $|\nu|=1$. Taking limits gives
$s_\nu(p)=0$, so $\nu\in N_{\mathcal K}(p)\subset V$.
Set
\begin{equation}\label{eq:v}
v=P_{V^\perp}(J\nu).
\end{equation}
Total reality of $V$ implies $J\nu\notin V$, hence $v\ne0$.
Moreover,
\[
\ip{\nu}{v}=0,
\qquad
\ip{J\nu}{v}=|v|^2,
\qquad
\ell_\nu(v)=i|v|^2.
\]

Kobayashi invariance and \eqref{eq:depth-harnack} give
\[
\begin{aligned}
k_U(z;Z_j(z))
&=k_U(\gamma_j(z);\delta_jv)\\
&\le\frac{\delta_j|v|}{\delta(\gamma_j(z))}
\le e^{2d_U^{\mathrm K}(o,z)}|v|.
\end{aligned}
\]
Together with the Euclidean lower bound in
\eqref{eq:kob-bounds}, this makes $Z_j$ locally uniformly
bounded. 

To prove that the limit has no zeros, introduce
\[
H_j(z)=1+
\frac{\ell_{\nu_j}(p_j-\gamma_j(z))}{\delta_j}.
\]
The normalization $s_{\nu_j}(p_j)=\delta_j$ and the chain
rule yield
\begin{equation}\label{eq:identities}
\operatorname{Re}H_j(z)
=\frac{s_{\nu_j}(\gamma_j(z))}{\delta_j}>0,
\qquad
H_j(o)=1,
\qquad
dH_j(Z_j)=-\ell_{\nu_j}(v).
\end{equation}
The Cayley transform of $H_j$ maps $U$ into $\D$ and
vanishes at $o$. Schwarz--Pick therefore gives
\[
\left|\frac{H_j(z)-1}{H_j(z)+1}\right|
\le\tanh d_U^{\mathrm K}(o,z),
\]
so $H_j$ is also locally uniformly bounded. Choose a common subsequence with $Z_j\to Z$ and $H_j\to H$
locally uniformly. Derivatives converge as well, and
\eqref{eq:identities} gives
\begin{equation}\label{eq:nonzero-flow}
\operatorname{Re}H\ge0,
\qquad
H(o)=1,
\qquad
dH(Z)=-i|v|^2.
\end{equation}
Since $\operatorname{Re}H$ is nonnegative harmonic and
$\operatorname{Re}H(o)=1$, the strong minimum principle gives
$\operatorname{Re}H>0$ throughout $U$.
Since $v\ne0$, the last identity excludes a zero of $Z$.
\end{proof}

\begin{lemma}[Completeness of the limiting field]
\label{lem:flow-complete}
The vector field $Z$ obtained in
Lemma~\ref{lem:flow-limit} has a complete real flow.
\end{lemma}

\begin{proof}
Keep the subsequence and the limiting field $Z$ fixed.
Let $z_*\in U$ and $T>0$, and use $x_j$ and $c$ as in
Lemma~\ref{lem:flow-segments}. That lemma ensures that
\[
y_j(t)=x_j+t\delta_jv
\]
belongs to $U$ for $|t|\le T$ and all sufficiently large $j$.
Hence the curves
\[
\psi_j(t)=\gamma_j^{-1}(y_j(t)),
\qquad |t|\le T,
\]
are well-defined. They satisfy
\[
 \psi_j(0)=z_*,
 \qquad
 \dot\psi_j(t)=Z_j(\psi_j(t)).
\]

The boundary-distance estimate gives
\[
k_U(y_j(t);\delta_jv)
\le\frac{\delta_j|v|}{\delta(y_j(t))}
\le\frac{2|v|}{c}.
\]
Integrating along the segment and using Kobayashi invariance,
we obtain
\begin{equation}\label{eq:trajectory}
d_U^{\mathrm K}(z_*,\psi_j(t))
=d_U^{\mathrm K}(x_j,y_j(t))
\le\frac{2T|v|}{c}.
\end{equation}
Put $C_T=2T|v|/c$. By Lemma~\ref{lem:gaps} and continuity of
$d_U^{\mathrm K}$, the closed ball
\[
 B_T=\overline{B}^{\mathrm K}_U(z_*,C_T)
   =\{z\in U:d_U^{\mathrm K}(z_*,z)\le C_T\}
\]
is closed in $U$ and lies in the relatively compact open ball
$B_U^{\mathrm K}(z_*,C_T+1)$; hence it is compact.
Equation~\eqref{eq:trajectory} places
every $\psi_j([-T,T])$ in this same set.

\begin{figure}[ht]
\centering
\setlength{\unitlength}{0.72mm}
\begin{picture}(150,39)
\thicklines
\qbezier(2,30)(29,38)(62,30)
\put(8,34){\makebox(0,0)[l]{$\partial U$}}
\put(4,7){\makebox(0,0){$U$}}
\put(13,17){\line(1,0){42}}
\put(32,17){\circle*{2.2}}
\put(32,12){\makebox(0,0){$x_j$}}
\put(39,13){\makebox(0,0)[l]{$y_j(t)$}}
\put(56,17){\vector(0,1){12}}
\put(54,25){\makebox(0,0)[r]{$\delta(y_j(t))\ge c\delta_j/2$}}
\put(66,21){\vector(1,0){17}}
\put(74,27){\makebox(0,0){$\gamma_j^{-1}$}}
\qbezier(85,20)(85,37)(116,37)
\qbezier(116,37)(147,37)(147,20)
\qbezier(147,20)(147,3)(116,3)
\qbezier(116,3)(85,3)(85,20)
\qbezier(95,20)(95,32)(116,32)
\qbezier(116,32)(137,32)(137,20)
\qbezier(137,20)(137,8)(116,8)
\qbezier(116,8)(95,8)(95,20)
\put(143,33){\makebox(0,0){$U$}}
\put(128,12){\makebox(0,0){$B_T$}}
\qbezier(98,18)(113,31)(133,19)
\put(112,25){\circle*{2.2}}
\put(112,29){\makebox(0,0){$z_*$}}
\put(121,24){\makebox(0,0)[l]{$\psi_j(t)$}}
\end{picture}
\caption{The supporting-gap estimate
keeps the translated segment at depth comparable to $\delta_j$.
After pullback by $\gamma_j^{-1}$, Kobayashi invariance keeps the
entire trajectory in the fixed compact ball $B_T\Subset U$.}
\label{fig:rescaled-trajectories}
\end{figure}

The fields $Z_j$ are uniformly bounded on $B_T$, so the curves
$\psi_j$ are uniformly Lipschitz. Arzel\`a--Ascoli gives
a subsequence converging uniformly on $[-T,T]$ to a curve
$\psi$. Uniform convergence of $Z_j$ on $B_T$ permits passage
to the limit in
\[
\psi_j(t)=z_*+\int_0^tZ_j(\psi_j(s))\,ds,
\]
giving
\[
\psi(t)=z_*+\int_0^tZ(\psi(s))\,ds.
\]
Therefore $Z$ has an integral curve through $z_*$ on
$[-T,T]$.

The subsequence of curves may depend on $z_*$ and $T$,
but the field $Z$ does not. Uniqueness of solutions makes
these integral curves agree on overlapping intervals.
Since $z_*$ and $T$ were arbitrary, the real flow of $Z$
is complete.
\end{proof}

\begin{proof}[Proof of Proposition~\ref{prop:flow}]
Lemma~\ref{lem:flow-limit} supplies a nowhere-zero
holomorphic vector field $Z$, and
Lemma~\ref{lem:flow-complete} proves its completeness.
Its flow maps $\Phi_t:U\to U$ are holomorphic by
holomorphic dependence on the initial point. Uniqueness
gives $\Phi_t^{-1}=\Phi_{-t}$, so
\[
 \Phi_t\in\Aut(U)
\]
for every $t\in\R$.
\end{proof}


\subsection{The convex factor}

We apply the complete flow first to $U$ and then to the complementary
factor in Theorem~\ref{thm:hhr}.

\begin{corollary}\label{cor:positive-aut}
A positive-dimensional bounded convex domain admitting a quotient
of finite K\"ahler--Einstein volume satisfies
\[
 \dim_\R\Aut(U)>0.
\]
\end{corollary}

\begin{proof}
Choose the complex extreme point with totally real normal span
provided by Proposition~\ref{prop:normal-point}. The global peak
theorem and Lemma~\ref{lem:peak} put it in an orbit limit set.
Proposition~\ref{prop:flow} then gives a nontrivial one-parameter
subgroup of $\Aut(U)$.
\end{proof}

\begin{rem*}
The argument also applies to compact quotients. Zimmer~\cite{Zim21} proved this corollary under the additional assumption that the convex domain is $C^{1,\alpha}$-smooth. Corollary~\ref{cor:positive-aut} requires no
boundary regularity and allows finite stabilizers.
\end{rem*}

To apply this corollary to the complementary factor in
Theorem~\ref{thm:hhr}, we recall Frankel's convex-factor argument.

\begin{lemma}\label{lem:convex-factors}
If $U\subset\C^n$ is bounded and convex and
$\Phi:X_1\times X_2\to U$ is biholomorphic, then the images
$\Phi(X_1\times\{x_2\})$ and $\Phi(\{x_1\}\times X_2)$ are bounded convex
domains in complex affine subspaces of $\C^n$.
\end{lemma}

\begin{proof}
The assertion is essentially Frankel's convex-factor
lemma~\cite[Corollary~5.5]{Fra89}. We use the following parameter
version of Cartan's iteration argument~\cite{Car32}.
If $D$ is a bounded domain and $H:P\times D\to D$ is a
holomorphic family, with $P$ connected and $H_{t_0}=\id_D$,
then $H_t=\id_D$ for every $t\in P$.
First take $P$ to be a disk centered at $t_0=0$. If the germ of
$H$ is nonconstant, let $k$ be the least index for which
\[
H_t(w)=w+t^k v(w)+O(t^{k+1})
\]
has $v\not\equiv0$. Iteration gives
\[
H_t^{\circ m}(w)=\underbrace{H_t \circ H_t \circ \dots \circ H_t}_{m}(w)=w+m t^k v(w)+O(t^{k+1}).
\]
In particular, the coefficient of $t^k$ is exactly $m v(w)$.
Fix a parameter circle contained in the disk. Since $D$ is bounded,
the coordinate functions of $H_t^{\circ m}(w)$ are bounded there
independently of $m$. Cauchy's estimate therefore bounds $m v(w)$
independently of $m$, a contradiction. For a
several-dimensional parameter, restrict to a complex line on which
the first nonzero homogeneous Taylor term does not vanish. Thus
$H_t=\id_D$ near $t_0$, and the identity theorem extends the
conclusion over the connected parameter manifold $P$.

For $a,b\in X_2$, define
\[
F_{a,b}(x_1)
=
\pi_{X_1}\Phi^{-1}
\left(\frac{\Phi(x_1,a)+\Phi(x_1,b)}2\right),
\]
where $\pi_{X_1}:X_1\times X_2\to X_1$ is the projection.
Convexity of $U$ makes this map well-defined. Its value is
the first product coordinate of the midpoint of
$\Phi(x_1,a)$ and $\Phi(x_1,b)$.
To apply the preceding observation on the bounded domain $U$,
consider
\[
\widetilde F_{a,b}\bigl(\Phi(x_1,x_2)\bigr)
=
\Phi\bigl(F_{a,b}(x_1),x_2\bigr).
\]
This is a holomorphic self-map of $U$, and it depends
holomorphically on $(a,b)\in X_2\times X_2$.
When $a=b$, we have $F_{a,a}=\id_{X_1}$, so
$\widetilde F_{a,a}=\id_U$. Since $X_2\times X_2$ is connected,
the parameter rigidity gives $\widetilde F_{a,b}=\id_U$
and therefore $F_{a,b}=\id_{X_1}$.

Thus the midpoint of any two points in
$\Sigma_{x_1}:=\Phi(\{x_1\}\times X_2)$ remains in $\Sigma_{x_1}$.
This slice is closed relative to $U$, being a fiber of
$\pi_{X_1}\circ\Phi^{-1}$. Repeated midpoint formation and
dyadic approximation along segments in $U$ therefore show
that $\Sigma_{x_1}$ is convex. Interchanging the factors proves
the same for $\Phi(X_1\times\{x_2\})$.

Finally, let $\Sigma$ be either slice and fix $p\in\Sigma$.
For every $q\in\Sigma$, the segment from $p$ to $q$ lies in
$\Sigma$, so its initial tangent vector $q-p$ belongs to
$T_p\Sigma$. Hence, with $A=p+T_p\Sigma$,
\[
\Sigma\subset A.
\]
For every $q\in\Sigma$, the inclusion $\Sigma\subset A$ gives
$T_q\Sigma\subset T_p\Sigma$. These spaces have the same complex
dimension, so equality holds. The inverse function theorem now shows
that $\Sigma$ is open in the complex affine space $A$. It is
connected, convex, and bounded, as required.
\end{proof}

\begin{proof}[Proof of Theorem~\ref{thm:convex}]
For $n=1$, the result follows from uniformization. Assume $n\ge2$.
Corollary~\ref{cor:positive-aut} gives $G\ne\{1\}$.
A bounded convex domain is contractible and HHR/USq.
Theorem~\ref{thm:hhr} gives $U\cong S\times Y$ with $\Aut(Y)$ discrete.
Proposition~\ref{prop:factor-metric-volume} gives finite-volume quotients of
both factors and complete K\"ahler--Einstein factor metrics.
Suppose that $\dim_\C Y>0$ and fix $s_0\in S$.
Lemma~\ref{lem:convex-factors} identifies $Y$
biholomorphically with the bounded convex domain
\[
 \Sigma=\Phi(\{s_0\}\times Y)
\]
in a complex affine subspace. Under this identification, $\Lambda$
is conjugated to a discrete subgroup of $\Aut(\Sigma)$, and
Proposition~\ref{prop:factor-metric-volume} gives
$\Vol(\Lambda\quot\Sigma,g_Y)<\infty$. The transported metric $g_Y$
is complete and satisfies $\Ric(g_Y)=-g_Y$; by uniqueness on the
bounded convex domain $\Sigma$, it is the normalized complete
K\"ahler--Einstein metric. Corollary~\ref{cor:positive-aut} gives
$\dim_\R\Aut(Y)>0$, contradicting the discreteness of $\Aut(Y)$.
Thus $Y$ is a point and $U\cong S$.
\end{proof}

\section{Bers--Griffiths domains and proof of Theorem~\ref{mt-3}}\label{sec:griffiths}

We establish a uniform squeezing estimate for universal covers
of curve families and use it to derive the automorphism conclusion
from Theorem~\ref{thm:hhr}. For the classical discreteness results
and their formulation for algebraic families, see
\cite[\S3]{GR18}. For the algebraic families in Griffiths'
construction, finite K\"ahler--Einstein volume and
Theorem~\ref{thm:hhr} then give Theorem~\ref{mt-3}.

\subsection{Curve families and the squeezing estimate}
\label{subsec:griffiths-families}

For a complex surface $X$, define
\begin{equation}\label{eq:gr-squeezing}
\sq_X(x)=\sup\bigl\{r>0:
 f\colon X\hookrightarrow\B^2\text{ is holomorphic},\quad
 f(x)=0,\quad \B_r^2\subset f(X)\bigr\}.
\end{equation}
The supremum is zero if there is no such embedding. Thus $X$ is
HHR/USq precisely when $\inf\limits_{x\in X}\sq_X(x)>0$~\cite{LSY04, Yeu09}. The embeddings
in \eqref{eq:gr-squeezing} are defined on all of $X$ and may depend
on $x$.

A \emph{finite-type hyperbolic curve} is a compact Riemann surface of
genus $g$ with $m$ distinct points removed, where $2g-2+m>0$.
A \emph{smooth family} of these curves means a surjective holomorphic
submersion $\pi\colon M\to C$ between connected complex manifolds, with finite-type hyperbolic curves of one fixed type $(g,m)$ over a finite-type hyperbolic curve and the following local compactification property.
Each point of $C$ has a neighborhood $W$ over which $\pi$ extends
to a proper holomorphic submersion
\[
\overline M_W\longrightarrow W
\]
with compact curve fibers, such that
$\overline M_W\setminus\pi^{-1}(W)$ consists of $m$ disjoint
holomorphic sections. The case $m=0$ is allowed. Ehresmann's theorem
for the compactified family, with the trivializations chosen to
preserve the sections, makes $M\to C$ a locally trivial smooth
fiber bundle. This is the curve-family hypothesis used by Griffiths
\cite[\S3]{Gri71}.

When $3g-3+m>0$, Griffiths proved that the universal cover
$D=\widetilde M$ is a contractible bounded domain of
holomorphy~\cite[Lemma~6.2]{Gri71}; in particular, it is
pseudoconvex. The only hyperbolic type omitted by this numerical
condition is $(g,m)=(0,3)$. Its marked Teichm\"uller space is a point,
so after pullback to $\widetilde C=\D$ the marked family is
holomorphically trivial and $\widetilde M\cong\D\times\D$.
We prove the following uniform strengthening in all cases.
\begin{theorem}\label{thm:griffiths-squeezing}
Let $\pi\colon M\to C$ be a smooth family of finite-type hyperbolic
curves, where $C$ is itself a finite-type hyperbolic curve. Then its
universal cover $D=\widetilde M$ is  $\mathrm{HHR/USq}$. More precisely,
\be \label{eq:gr-main-bound}
\inf\limits_{p\in D}\sq_D(p)\ge1/5.
\ene
\end{theorem}

\begin{corollary}\label{cor:griffiths-zariski}
Let $X$ be a smooth quasi-projective complex surface and let $x\in X$.
There is a Zariski open neighborhood $M$ of $x$, supplied by
Griffiths' curve-family construction, whose universal cover satisfies
Theorem~\ref{thm:griffiths-squeezing}.
\end{corollary}

Here a \emph{Bers--Griffiths domain in complex dimension two}
means a bounded realization of the universal cover obtained from
that construction. Thus the theorem addresses the surface case of
\cite[Question~1.18]{Zim26}. The estimate uses the curve-family
structure; a finite-volume quasi-projective quotient alone is not
the hypothesis of the theorem. The proof combines a centered Bers realization with elementary
one-variable estimates. Only the base and the fibers are
one-dimensional; no bound on the Teichm\"uller-space dimension is used.

\subsection{A centered global Bers model}

Choosing the reference fiber at the point under consideration gives
a global model whose central fiber is the unit disk. Write $\D^*=\{\zeta\in\widehat{\C}:|\zeta|>1\}$, including $\infty$.

\begin{lemma}\label{lem:gr-normalization}
For each $p\in D$, there is a biholomorphism
\begin{equation}\label{eq:gr-model}
\Theta_p\colon D\longrightarrow
\Omega=\{(z,w):z\in\D,\ w\in \Omega_z\}\subset\D\times\C
\end{equation}
with $\Theta_p(p)=(0,0)$ and $\Omega_0=\D$.
Each $\Omega_z$ is a bounded Jordan domain with quasicircle boundary.
There are conformal maps
\begin{equation}\label{eq:gr-exterior}
f_z\colon\D^*\longrightarrow\widehat{\C}\setminus\overline{\Omega_z},
\qquad f_z(\zeta)=\zeta+O(\zeta^{-1}),\qquad f_0=\id,
\end{equation}
which depend holomorphically on $(z,\zeta)$ for $|\zeta|>1$
and extend homeomorphically to $|\zeta|=1$. 
\end{lemma}

\begin{proof}
Pull the family back by the universal covering $\widetilde C=\D\to C$,
choosing the base coordinate so that the image of $p$ lies over zero.
Write $M_{\D}=M\times_C\D$. On overlaps, the identity maps between
the punctured families extend uniquely across the missing points.
The local compactifications therefore glue to a proper holomorphic
submersion $\overline M_{\D}\to\D$. Its puncture locus is an
unramified $m$-sheeted covering of $\D$ and, because $\D$ is simply
connected, splits into $m$ disjoint global sections. Let
$\overline\Sigma_0$ be the compact central fiber, let $P_0$ be its
set of marked points, and put
\[
                 \Sigma_0=\overline\Sigma_0\setminus P_0.
\]
Thus $\Sigma_0$ is the central fiber of $M_{\D}$.

Ehresmann's theorem and the contractibility of $\D$ provide a
marking of the pair $(\overline\Sigma_0,P_0)$, preserving the
puncture labels.
Choose a torsion-free cofinite Fuchsian group
$\Delta<\Aut(\D)$ and a pointed uniformization
\[
                         \Sigma_0\cong\Delta\quot\D
\]
taking zero to the image of $p$. The universal property of
Teichm\"uller space gives a holomorphic classifying map
\[
\tau:\D\longrightarrow\mathcal T(\Delta),
\qquad \tau(0)=0.
\]
This is the empty-multicurve case of \cite[Theorem~10.1]{HK14};
in the present setting it is also the classifying map constructed
in \cite[Lemma~6.1]{Gri71}.

For $\tau(z)=[\mu]$, extend the $\Delta$-invariant Beltrami
coefficient $\mu$ by zero to $\D^*$. Let $W^\mu$ be its
principal quasiconformal solution, normalized by
\[
W^\mu(\zeta)=\zeta+O(\zeta^{-1})
\qquad (\zeta\to\infty),
\]
and set
\[
\Omega_z=W^\mu(\D),
\qquad
f_z=W^\mu|_{\D^*}.
\]
Since $\mu$ vanishes on $\D^*$, the map $f_z$ is conformal there.
The sphere homeomorphism $W^\mu$ shows that $\Omega_z$ is a bounded
quasidisk and supplies the boundary extension of $f_z$. These
definitions are independent of the representative $\mu$: equivalent
representatives give exterior maps differing by a M\"obius
transformation, while the normalization at infinity forces that
transformation to be the identity~\cite{Ber73}.

The principal normalization is compatible with holomorphic
dependence. On a sufficiently small disk in the base, choose a local
holomorphic section of the Teichm\"uller projection above $\tau$, and
denote its Beltrami coefficients by $\mu_z$. Let $\widehat W_z$ be the
corresponding Ahlfors--Bers solution normalized to fix $\infty$. Write
\[
 \widehat W_z(\zeta)=a(z)\zeta+b(z)+O(\zeta^{-1}),
 \qquad a(z)\ne0.
\]
The functions $a$ and $b$ are holomorphic, and
\[
                         W_z=\frac{\widehat W_z-b(z)}{a(z)}
\]
has the required principal normalization. Ahlfors--Bers
holomorphic dependence~\cite{AB60} now shows that
$(z,\zeta)\mapsto f_z(\zeta)$ is jointly holomorphic for
$|\zeta|>1$. Independence of the representative makes these local
constructions agree on overlaps.

Relative Bers uniformization now gives a holomorphic covering
\[
\Omega=\{(z,w):z\in\D,\ w\in\Omega_z\}
\longrightarrow M_{\D},
\]
whose restriction to each fiber is its universal covering;
see \cite{Ber73} and
\cite[Proposition~5.8 and Lemma~6.1]{Gri71}.
Local continuous lifts of $\tau$ give topological trivializations
\[
(z,v)\longmapsto\bigl(z,W^{\mu_z}(v)\bigr).
\]
Thus $\Omega\to\D$ is a topological disk bundle over the
contractible base $\D$, so $\Omega$ is contractible and, in
particular, simply connected. Consequently
$\Omega\to M_{\D}\to M$ is the universal covering of $M$.
At $z=0$ we may take $\mu=0$, yielding
$\Omega_0=\D$ and $f_0=\id$.
By our choice of the central uniformization, $(0,0)$ and $p$
lie over the same point of $M$. Uniqueness of the pointed
universal covering therefore gives the required biholomorphism
$\Theta_p:D\to\Omega$ with $\Theta_p(p)=(0,0)$.
\end{proof}

\subsection{Uniform squeezing}

The second-coefficient estimate for univalent functions bounds every
fiber, and Schwarz's lemma supplies a uniform inner neighborhood.

\begin{lemma}\label{lem:gr-inclusions}
The model $\Omega$ in Lemma~\ref{lem:gr-normalization} satisfies
\begin{equation}\label{eq:gr-inclusions}
\{(z,w):3|z|+|w|<1\}\subset \Omega\subset\D\times\D_2.
\end{equation}
\end{lemma}

\begin{proof}
\underline{Outer inclusion.} Fix $z\in\D$ and $w\in \Omega_z$. The function
\[
h(\eta)=\frac{1}{f_z(1/\eta)-w}
        =\eta+w\eta^2+O(\eta^3)
\]
extends holomorphically across zero and is univalent on $\D$.
Indeed, its denominator never vanishes, it is injective away from
zero, and zero is its only zero. The classical second-coefficient
estimate~\cite[Chapter~2]{Dur83} for normalized univalent functions
gives $|w|\le2$. Since this holds for every $w\in\Omega_z$ and
$\Omega_z$ is open, it follows that
$\Omega_z\subset\D_2$ and
$\overline{\Omega_z}\subset\overline{\D_2}$.

\smallskip
\underline{Inner inclusion.}
We first estimate the displacement of the exterior maps.
For fixed $z\in\D$, define
\[
q_z(\eta)=f_z(1/\eta)-\frac1\eta
\qquad (0<|\eta|<1).
\]
The normalization of $f_z$ shows that $q_z$ extends
holomorphically across zero with $q_z(0)=0$.
It is continuous on $\overline{\D}$, and on $|\eta|=1$
the boundary correspondence and the outer inclusion give
\[
|q_z(\eta)|
\le |f_z(1/\eta)|+1
\le 3.
\]
The maximum modulus principle therefore yields
\[
|f_z(\zeta)-\zeta|\le3
\qquad (z\in\D,\ |\zeta|\ge1).
\]

Now fix a finite $\zeta$ with $|\zeta|>1$ and consider
the base variable $z$. The function
\[
\varphi_\zeta(z)=\frac{f_z(\zeta)-\zeta}{3}
\]
is holomorphic on $\D$, has modulus at most $1$, and
vanishes at zero because $f_0=\id$.
Schwarz's lemma gives
\begin{equation}\label{eq:gr-exterior-schwarz}
|f_z(\zeta)-\zeta|\le3|z|.
\end{equation}
For each fixed $z$, continuity of $f_z$ up to the unit
circle extends this estimate to $|\zeta|=1$. Suppose that
$3|z|+|w|<1$.
If $w\notin\Omega_z$, the exterior map and its boundary
extension give a finite $\zeta$ with $|\zeta|\ge1$
such that $w=f_z(\zeta)$.
Indeed, exterior points are covered by $f_z$ on $\D^*$,
and boundary points by its extension to the unit circle.
Using \eqref{eq:gr-exterior-schwarz}, we obtain
\[
|w|
=|f_z(\zeta)|
\ge|\zeta|-|f_z(\zeta)-\zeta|
\ge1-3|z|,
\]
contrary to the choice of $(z,w)$.
Hence $w\in\Omega_z$, proving the inner inclusion.
\end{proof}

\begin{proof}[Proof of Theorem~\ref{thm:griffiths-squeezing}]
Fix $p\in D$. Lemma~\ref{lem:gr-normalization} gives a biholomorphism
$\Theta_p\colon D\to \Omega$ with $\Theta_p(p)=(0,0)$.
Set
\[
A(z,w)=\left(\sqrt{\frac35}\,z,\frac{w}{\sqrt{10}}\right),
\qquad \Psi_p=A\circ\Theta_p.
\]
The proof is summarized by the inclusion diagram
\[
\begin{array}{c}
\B_{1/5}^2
 \subset A\bigl(\{(z,w):3|z|+|w|<1\}\bigr)
 \subset A(\Omega) \\[3pt]
 A(\Omega)\subset A(\D\times\D_2)\subset\B^2.
\end{array}
\]
We verify the two outer inclusions. The map $\Psi_p$ is a
holomorphic embedding of the entire universal cover. For
$(z,w)\in \Omega$, Lemma~\ref{lem:gr-inclusions} gives
\[
\lVert A(z,w)\rVert^2
=\frac{6|z|^2+|w|^2}{10}<\frac{6+4}{10}=1.
\]
Hence $\Psi_p(D)\subset\B^2$. Conversely, if
$v=(v_1,v_2)\in\B_{1/5}^2$ and $(z,w)=A^{-1}(v)$, then
Cauchy--Schwarz gives
\[
3|z|+|w|
=\sqrt{15}\,|v_1|+\sqrt{10}\,|v_2|
\le5\lVert v\rVert<1.
\]
The lower inclusion in \eqref{eq:gr-inclusions} implies
$(z,w)\in \Omega$, so $v\in\Psi_p(D)$. Since $p\in D$ is arbitrary, this proves
\eqref{eq:gr-main-bound}.
\end{proof}

\begin{proof}[Proof of Corollary~\ref{cor:griffiths-zariski}]
Griffiths' construction provides a Zariski open neighborhood $M$
of $x$ and a smooth algebraic fibration $M\to C$ with connected
finite-type fibers~\cite[Lemma~2.2]{Gri71}. In complex dimension
two, the base $C$ is a Zariski open subset of the affine line arising
from the generic projection used in the proof of
\cite[Lemma~2.1]{Gri71}; in particular, it is noncompact. Moreover,
\cite[Lemma~2.2(iv)]{Gri71} realizes $\widetilde C$ as a bounded
simply connected domain in $\C$, and hence $\widetilde C\cong\D$.
The fibers have fixed genus and fixed number of punctures, while
\cite[Lemma~3.2]{Gri71} supplies the required local
compactifications. Theorem~\ref{thm:griffiths-squeezing} therefore
applies to $M$.
\end{proof}

\begin{rem*}
If the marked classifying map is constant, the centered model is
$\Omega=\D\times\D$, and $f_z=\id$ for every $z$. Thus the bidisk
fits the construction directly. In a varying family, the estimates
remain uniform even when the fibers degenerate at an end of the
base: they require no common bound on the quasiconformal dilatations.
All fibers over points of $C$ are smooth, as stipulated in the theorem. Corollary~\ref{cor:griffiths-zariski} concerns the universal cover of the
Griffiths neighborhood $M$. It does not identify that cover with
the universal cover of the original surface $X$.
\end{rem*}

\subsection{Automorphism rigidity}
\label{subsec:griffiths-rigidity}

\begin{proof}[Proof of Theorem~\ref{mt-3}]
Let $\pi:D\to M$ be the universal covering associated with a
Griffiths family, and suppose that $D$ is not biholomorphic to
$\D\times\D$. By Theorem~\ref{thm:griffiths-squeezing}, $D$ is a
contractible bounded HHR/USq domain. Since it is pseudoconvex, it
carries a unique complete K\"ahler--Einstein metric
$g_D$~\cite{MY83}, normalized by
\[
                              \Ric(g_D)=-g_D.
\]
Uniqueness makes every deck transformation an isometry
of $g_D$; hence $g_D$ descends to a K\"ahler metric $g_M$ on $M$
with the same Ricci identity. Since $M$ is smooth and
quasi-projective, Griffiths' volume theorem
\cite[Proposition~7.3]{Gri71} gives
\[
                              \Vol_{g_M}(M)<\infty.
\]
Thus all hypotheses of Theorem~\ref{thm:hhr} are satisfied, and
\[
                              D\cong S\times Y,
\]
where $S$ is bounded symmetric and $Y$ is contractible with
$\Aut(Y)$ discrete.

Assume that $\dim_\C S>0$. Since $\dim_\C D=2$, either
\[
 (\dim_\C S,\dim_\C Y)=(1,1)
 \qquad\text{or}\qquad
 (\dim_\C S,\dim_\C Y)=(2,0).
\]
In the first case, $Y$ is a contractible Riemann surface.
Uniformization identifies it with either $\C$ or $\D$; the sphere
is excluded by contractibility. Both have nondiscrete automorphism
groups, contradicting the defining property of $Y$.

Hence $Y$ is a point and $\dim_\C S=2$. The classification of
two-dimensional bounded symmetric domains leaves only $\B^2$ and
$\D\times\D$. The base curve of a Griffiths neighborhood is
noncompact, so \cite[Theorem~1]{IN05} excludes $\B^2$ as the
universal cover, while the bidisk is excluded by hypothesis. This
contradiction shows that $S$ is a point. Therefore $D\cong Y$, and
$\Aut(D)$ is discrete.
\end{proof}


\end{document}